\documentclass{article}
\begin{document}

\centerline{\textbf{On certain theta functions and modular forms}} 
\centerline{\textbf{in Ramanujan theories}}
\[
\]
\centerline{\bf N.D. Bagis}
\centerline{\bf Aristotele University of Thessaloniki}
\centerline{\bf Thessaloniki Greece}
\centerline{\bf nikosbagis@hotmail.gr}
\[
\]
\begin{quote}

\centerline{\bf abstract\rm}
We give evaluations of certain Borwein's theta functions which appear in Ramanujan theory of alternative elliptic modular bases. Most of this theory where developed by B.C. Berndt, S. Bhargava and F.G. Garvan. We also study the most general class of these theta functions and give evaluation theorems and conjectures.   

\[\]     

\bf keywords: \rm{Ramanujan; theta functions; alternative modular bases; Evaluations;}

\end{quote}

\section{Introduction}

The Borwein-Borwein theorem for the cubic modular base ${}_2F_1\left(\frac{1}{3},\frac{2}{3};1;x\right)$, states that (see [6]):\\
\\
\textbf{Theorem 1.}\\ 
Let $\omega=e^{2\pi i/3}$, $q=e^{-\pi\sqrt{r}}$, $r>0$ and
\begin{equation}
a_0(q):=\sum^{\infty}_{m,n=-\infty}q^{m^2+mn+n^2},
\end{equation}
\begin{equation}
b_0(q):=\sum^{\infty}_{m,n=-\infty}\omega^{m-n}q^{m^2+mn+n^2},
\end{equation}
and
\begin{equation}
c_0(q):=\sum^{\infty}_{m,n=-\infty}q^{(m+1/3)^2+(m+1/3)(n+1/3)+(n+1/3)^2}.
\end{equation}
Then the solution of equation
\begin{equation}
\frac{{}_2F_1\left(\frac{1}{3},\frac{2}{3};1,1-x\right)}{{}_2F_1\left(\frac{1}{3},\frac{2}{3};1,x\right)}=\sqrt{r}
\end{equation}
is
\begin{equation}
x=\alpha_{r}=\left(\frac{c_0\left(q^{2/\sqrt{3}}\right)}{a_0\left(q^{2/\sqrt{3}}\right)}\right)^3
\end{equation}
and
\begin{equation}
a_0^3(q)=b^3_0(q)+c_0^3(q)
\end{equation}
\\
Also in view of [6] which the case of Ramanujan elliptic modular bases theory is treated extensively, many results can be extracted. In this work we explore some of these ideas. From our part we give formulas for closed form evaluations of the above theorem theta functions.

\section{Some Results}

It is known from [6] p.82, the following relation of hypergeometric functions  
\begin{equation}
\sqrt[4]{1+8x}\cdot{}_2F_1\left(\frac{1}{3},\frac{2}{3};1,x\right)={}_2F_1\left(\frac{1}{6},\frac{5}{6};1;\frac{1}{2}-\frac{1-20x-8x^2}{2(1+8x)^{3/2}}\right).
\end{equation}
Setting $x=\alpha_{3r}$ in (7) and using (see [2] Theorem 4)
\begin{equation}
\beta_r=\frac{1}{2}-\frac{1-20\alpha_{3r}-8\alpha_{3r}^2}{2(1+8\alpha_{3r})^{3/2}},
\end{equation}
we get\\
\\
\textbf{Theorem 1.}\\
Let $r$ be real positive, then
\begin{equation}
u_r={}_2F_{1}\left(\frac{1}{6},\frac{5}{6};1;\beta_r\right)=\sqrt[4]{1+8\alpha_{3r}}\cdot{}_2F_{1}\left(\frac{1}{3},\frac{2}{3};1;\alpha_{3r}\right)
\end{equation}

Using the identities $\alpha_{1/r}=\alpha'_r$ and $\beta_{1/r}=\beta'_r$ one can easily see that the multiplier of 3rd-degree in the 6-base (signature 6) is\\ 
\\
\textbf{Theorem 2.}
\begin{equation}
m_3(6)=\frac{{}_2F_{1}\left(\frac{1}{6},\frac{5}{6};1;\beta_{9r}\right)}{{}_2F_{1}\left(\frac{1}{6},\frac{5}{6};1;\beta_r\right)}=\sqrt{3r}\cdot\sqrt[4]{\frac{1+8\alpha'_{3r}}{1+8\alpha_{3r}}}
\end{equation} 

From [6] Theorem 2.10 we have\\
\\
\textbf{Theorem 3.}\\If $q=e^{-\pi\sqrt{r}}$, $r>0$ and $z_r={}_2F_1\left(\frac{1}{3},\frac{2}{3};1;\alpha_{r}\right)$, then
\begin{equation}
z_{3r}=a_0(q^2).
\end{equation}
\\
\textbf{Lemma 1.}\\
If $x>0$, then
\begin{equation}
\sum^{\infty}_{n=1}\frac{X(n)}{e^{nx}-1}=\sum^{\infty}_{n=1}e^{-nx}\sum_{d|n}X(d),
\end{equation}
when both sides of (12) are convergent.\\
\\
\textbf{Proof.}\\
$$
\sum^{\infty}_{n=1}\frac{X(n)}{e^{nx}-1}=\sum^{\infty}_{n=1}\frac{X(n)e^{-nx}}{1-e^{-nx}}=\sum^{\infty}_{n=1}X(n)\sum^{\infty}_{m=1}e^{-mx}=
$$
$$
=\sum^{\infty}_{n,m=1}e^{-nmx}X(n)=\sum^{\infty}_{n=1}e^{-nx}\sum_{d|n}X(d).
$$
\\
\textbf{Corollary 1.}\\The number of representations $s(n)$, of the integer $n$ in the form
\begin{equation}
n=x^2+y^2+xy+zw+z^2+w^2
\end{equation}
is
\begin{equation}
s(n)=12\sum_{d|n}d-12\sum_{d|3,d|n}d
\end{equation}
and $s(0)=1$.\\
\\
\textbf{Proof.}\\
In [5], p.460, Entry 3(i) it has shown that 
\begin{equation}
a_0(q)^2=1+12\sum^{\infty}_{n=1}\frac{nq^n}{1-q^n}-36\sum^{\infty}_{n=1}\frac{nq^{3n}}{1-q^{3n}} 
\end{equation}
Using (1) along with the formula (Lemma 1)
\begin{equation}
\sum^{\infty}_{n=1}\frac{X(n)q^n}{1-q^n}=\sum^{\infty}_{n=1}\left(\sum_{d|n}X(d)\right)q^n,
\end{equation} 
we get the result.\\
\\
\textbf{Theorem 4.}\\ The number of representations of the integer $n$ in the form
\begin{equation}
n=x^2+xy+y^2
\end{equation}
is
\begin{equation}
r(n)=6\sum_{d\equiv1(3),d|n}1-6\sum_{d\equiv2(3),d|n}1
\end{equation}
and $r(0)=1$.\\
\\
\textbf{Proof.}\\ 
From the proof of [6] Theorem 2.12 p.11-12 and (11),(1) of the present article, we have 
\begin{equation}
\sum^{\infty}_{m,n=-\infty}q^{n^2+mn+m^2}=1+6\sum^{\infty}_{n=0}\frac{q^{3n+1}}{1-q^{3n+1}}-6\sum^{\infty}_{n=0}\frac{q^{3n+2}}{1-q^{3n+2}}.
\end{equation}
Expanding the above formula into power series using (12) we get the result.\\

For the 5th base we have the relation\\
\\
\textbf{Theorem 5.}\\If $q=e^{-\pi\sqrt{r}}$, $r>0$ we have
\begin{equation}
u_r^4=u(\beta_r)^4=1+240\sum^{\infty}_{n=0}\frac{n^3q^{2n}}{1-q^{2n}}=\frac{1}{2}\left(a^8+b^8+c^8\right)
\end{equation}
where 
$a=\theta_2(q)=\sum^{\infty}_{n=-\infty}q^{(n+1/2)^2}$, $b=\theta_3(q)=\sum^{\infty}_{n=-\infty}q^{n^2}$ and  $c=\theta_4(q)=\sum^{\infty}_{n=-\infty}(-1)^nq^{n^2}$, $|q|<1$.\\

From (15) and [3],[4] we also have
$$
a_0(q)^2=1+12\sum^{\infty}_{n=1}\frac{nq^n}{1-q^n}-36\sum^{\infty}_{n=1}\frac{nq^{3n}}{1-q^{3n}}=
$$
\begin{equation} 
=1-12q\frac{d}{dq}\log\left(\eta(q^3)^{-1}\vartheta\left(\frac{3}{2},\frac{1}{2};q\right)\right) 
\end{equation}
where
\begin{equation}
\vartheta\left(\frac{3}{2},\frac{1}{2};q\right)=q^{1/12}\eta(q^3)Q(k_r^2)^{1/12}
\end{equation}
The minimal polynomial of $u=Q(v)$ is  
$$
-16 u^3 + 387420489 v + 19131876 u v + 196830 u^2 v + 84 u^3 v +u^4 v-1549681956 v^2-
$$
$$
-76527504 u v^2-787320 u^2 v^2-12480 u^3 v^2-4u^4 v^2+2324522934 v^3+114791256 u v^3+
$$
$$
+1180980 u^2 v^3- 40712 u^3 v^3+6 u^4 v^3-1549681956 v^4-76527504 u v^4-787320 u^2 v^4-
$$
$$
-12480 u^3 v^4-4 u^4 v^4+387420489 v^5+19131876 u v^5+196830 u^2 v^5+
$$
\begin{equation}
+84 u^3 v^5+u^4 v^5-16 u^3 v^6=0
\end{equation}
and is solvable. Hence 
\begin{equation}
a_0(q)^2=-q\frac{d}{dq}\log Q(k_r^2)
\end{equation}
hence
\begin{equation}
a_0(q)^2=-q\frac{d}{dk}\left(\log Q(k^2)\right)\frac{dk}{dq}
\end{equation}
But from [5] chapter 17 Entry 9 p.120 we have
\begin{equation}
\frac{dq}{dk}=\frac{q\pi^2}{2k(1-k^2)K(k)^2}.
\end{equation}
Setting (here the prime on $Q$ means derivative): 
\begin{equation}
\lambda_{r}:=\sqrt{-\frac{Q'(k_r^2)}{Q(k_r^2)}}
\end{equation}
we get the next\\
\\
\textbf{Theorem 6.}\\ 
If $u=Q(v)$ is that of (23) and $q=e^{-\pi\sqrt{r}}$, then
\begin{equation}
a_0(q)=\sum^{\infty}_{n,m=-\infty}q^{m^2+mn+n^2}=\frac{2K}{\pi}k_rk'_r\lambda_r
\end{equation}
Also
\begin{equation}
a_0(q^2)=\sum^{\infty}_{n,m=-\infty}q^{2[m^2+mn+n^2]}=\frac{2K}{\pi}\sqrt[4]{\frac{1-k_r^2+k_r^4}{1+8\alpha_{3r}}}
\end{equation}
Moreover if $r-$positive rational, then
\begin{equation}
\frac{\pi a_0(q)}{K}=Algebraic
\end{equation}
where $k'=\sqrt{1-k^2}$ is the complementary modulus.\\
\\
\textbf{Example 1.}\\ 
Let $q=e^{-\pi}$, then $k_1=\frac{1}{\sqrt{2}}$. Solving equation (23) (which is quartic) and setting $v=1/2$ we get 
$$
Q\left(\frac{1}{2}\right)=81 \left(885 + 511 \sqrt{3} - 3 \sqrt{174033 + 100478 \sqrt{3}}\right) 
$$
and
$$
Q'\left(\frac{1}{2}\right)=162 \left(5082 + 2934 \sqrt{3} - \sqrt{51655599 + 29823374 \sqrt{3}}\right)
$$
Also
$$
K(k_1)=K\left(\frac{1}{\sqrt{2}}\right)=\frac{8\pi^{3/2}}{\Gamma\left(-\frac{1}{4}\right)^2}
$$
Hence
\begin{equation}
a_0\left(e^{-\pi}\right)=\frac{8 \sqrt{2  \pi} \sqrt{\sqrt{24+14 \sqrt{3}}-3}}{\sqrt[4]{3} \Gamma \left(-\frac{1}{4}\right)^2}
\end{equation}
\\

From relations (5),(6) and (28) we get the next\\
\\
\textbf{Theorem 7.}\\ 
If $q=e^{-\pi \sqrt{r}}$ with $r>0$, then 
$$
c_0(q^2)=\sum^{\infty}_{m,n=-\infty}q^{2[(m+1/3)^2+(m+1/3)(n+1/3)+(n+1/3)^2]}=
$$
\begin{equation}
=\frac{2K}{\pi}\frac{1-k'_r}{1+k'_r}\sqrt{k'_r}\alpha_{3r}^{1/3}\lambda_{4r}
\end{equation}
Also
\begin{equation}
b_0(q^2)=\sum^{\infty}_{m,n=-\infty}\omega^{m-n}q^{2[m^2+mn+n^2]}=\frac{2K}{\pi}\frac{1-k'_r}{1+k'_r}\sqrt{k'_r}\lambda_{4r}\sqrt[3]{1-\alpha_{3r}}
\end{equation}
and consequently\\
\\
\textbf{Theorem 8.}\\ 
For $r>0$, we have
\begin{equation}
\left(1+8\alpha_{3r}\right)\lambda_{4r}^4=\frac{({k'}_r^2+k^{4}_r)\left(1+k'_r\right)^4}{k^{{'}2}_r(1-k'_r)^4}
\end{equation}
\\
\textbf{Theorem 9.}\\ The equation
\begin{equation}
t^4+8C_1t^3+18t^2-27=0
\end{equation}
have solution
\begin{equation}
t=\sqrt{1+8\alpha_{3r}}
\end{equation}
where $r$ is given from $C_1=1-2\beta_r$.\\
\\
\textbf{Proof.}\\ 
From (8) we get that equation
\begin{equation}
\frac{1-20x-8x^2}{(1+8x)^{3/2}}=C_1
\end{equation}
admits solution $x=\alpha_{3r}$ with $r$ such that $C=1-2\beta_r$.\\
\\
\textbf{Note 2.}\\ 
\textbf{i)} It holds
\begin{equation}
j_r=\frac{432}{\beta_r(1-\beta_r)}\textrm{, }r>0
\end{equation}
where $j_r$ is the $j-$invariant. Also
\begin{equation}
j_{r}=256\frac{\left((k'_r)^2+k_r^4\right)^3}{(k_rk'_r)^4}.
\end{equation}
Hence we have evaluated $\beta_r$, $\alpha_r$, $\lambda_r$ (see [2] and [7]) in terms of $k_r$.

\section{Binary theta functions}

\textbf{Definition 1.}\\ 
Assume the quadratic forms 
\begin{equation}
Q_1=a_1x^2+b_1xy+c_1y^2\textrm{ and }Q_2=a_2x^2+b_2xy+c_2y^2,
\end{equation}
where $a_1,b_1,c_1,a_2,b_2,c_2$ integers. The two forms $Q_1$ and $Q_2$ are equivalent iff exist   
$$
g=\left(
\begin{array}{cc}
	r\textrm{ }s\\
	t\textrm{ }u
\end{array}
\right)\in GL_2\left(\textbf{Z}\right)
$$
such that $Q_2=gQ_1$, where 
\begin{equation}
gQ_1:=a_1(rx+ty)^2+b_1(rx+ty)(sx+uy)+c_1(sx+uy)^2
\end{equation} 
\\
\textbf{Note 3.}\\Actually it is $g\in GL_2\left(\textbf{Z}\right)$ iff $r,t,s,u\in\textbf{Z}$ and $det(g)=\pm 1$.\\
\\
\textbf{Definition 2.}\\
Let $a_1,b_1,c_1,a_2,b_2,c_2$ be integers with $a_1,a_2,c_1,c_2>0$ and $D=b_1^2-4a_1c_1=b_2^2-4a_2c_2=-|D|<0$. 
If also 
\begin{equation}
\sum^{\infty}_{n,m=-\infty}q^{a_1n^2+b_1nm+c_1m^2}=\sum^{\infty}_{n,m=-\infty}q^{a_2n^2+b_2nm+c_2m^2}\textrm{, }\forall q: |q|<1,
\end{equation}  
then $Q_1$ and $Q_2$ are called $\theta-$equivalent and we write 
\begin{equation}
(a_2,b_2,c_2)\stackrel{\theta}{\equiv}(a_1,b_1,c_1)
\end{equation}
\\
\textbf{Proposition 1.}\\
If $Q_1$ and $Q_2$ are equivalent with $a_1,a_2,c_1,c_2>0$ and  $b_1^2-4a_1c_1,b_2^2-4a_2c_2<0$, then $Q_1$, $Q_2$ are also $\theta-$equivalent.\\
\\
\textbf{Theorem 10.}\\ 
Let $a,c>0$ and
\begin{equation}
\phi(q)=\phi(a,b,c;q):=\sum^{\infty}_{n,m=-\infty}q^{an^2+bnm+cm^2}\textrm{, }q=e^{2\pi i z}\textrm{, }Im(z)>0.
\end{equation}
Then for all $\theta-$equivalent functions $\phi(q)$ with $D=b^2-4ac=-|D|\in\left\{-3,-4,-7\right\}$, we have 
\begin{equation}
\phi(q)^2=\frac{1}{D+1}\left(1-24 \sum^{\infty}_{n=1}\frac{nq^n}{1-q^n}+D\left(1-24\sum^{\infty}_{n=1}\frac{nq^{-Dn}}{1-q^{-D n}}\right)\right).
\end{equation} 
\\
\textbf{Proof.}\\
The case $D=-3$ is relation (15) (see related references). For the cases $D=-4$ and $D=-7$, see Theorems 12,13 and the examples and notes below.\\
\\

From the above Theorem we get the following\\
\\
\textbf{Corollary 2.}\\  
The number of representations of the integer $n\neq0$ in the form
\begin{equation}
a(x^2+z^2)+b(xy+zw)+c(y^2+w^2),
\end{equation}
with $D=b^2-4ac=-|D|=-3,-4,-7$, is
\begin{equation}
s(n)=-\frac{24}{D+1}\sum_{d|n}d+\frac{24}{D+1}\sum_{d|n,-D|d}d.
\end{equation}
\\

In the case of $D=-3$, we set first (in general)
\begin{equation}
k_i(x):=\left(\frac{K\left(\sqrt{1-x^2}\right)}{K(x)}\right)^2.
\end{equation}
For all $(a,b,c)\stackrel{\theta}{\equiv}(1,1,1)$ obviously we have 
\begin{equation}
\sum^{\infty}_{n,m=-\infty}q^{an^2+bnm+cm^2}={}_2F_1\left(\frac{1}{3},\frac{2}{3};1;\alpha_{3r/4}\right).
\end{equation}
Also for all $r=k_i\left(\frac{m}{n}\right)$, $m<n$, $m,n$ positive integers, the minimal polynomial of
\begin{equation}
S(x)=\left(\frac{\pi}{2K(x)}\sum^{\infty}_{n,m=-\infty}q_1^{an^2+bnm+cm^2}\right)^2\textrm{, }(a,b,c)\stackrel{\theta}{\equiv}(1,1,1)
\end{equation}
with $q_1=e^{-\pi\sqrt{k_i(x)}}$, $x\in\textbf{Q}$, $0<x<1$,
has always degree equivalent to 4 and hence is solvable, since from Theorem 6 and (27),(23) it is
\begin{equation}
S(x)=-x(1-x^2)\frac{Q'\left(x^2\right)}{Q\left(x^2\right)}. 
\end{equation}

Moreover set in general
\begin{equation}
S(x):=\left(\frac{\pi}{2K(x)}\sum^{\infty}_{n,m=-\infty}q_1^{an^2+bnm+cm^2}\right)^2\textrm{, }q_1=e^{-\pi\sqrt{k_i(x)}}\textrm{, }0<x<1
\end{equation}
then\\
\\
\textbf{Theorem 11.}(Conjecture)\\For given integers $a,b,c$, with $a,c>0$ and $b^2-4ac<0$, the degree of $S(x)$ is constant when $x\in\textbf{Q}^{*}_+$, $0<x<1$.\\One can experimentally see that always exist integer coefficients $\alpha^{*}_{mn}$ and $N\in\textbf{N}$ such that
\begin{equation}
\sum^{N}_{n,m=0}\alpha^{*}_{nm}\left(S(x)\right)^nx^m=0\textrm{, }x\in\textbf{R}\textrm{, }0<x<1,
\end{equation}
\\
\textbf{Note 4.}\\ The set $P=\left\{x:x\in\textbf{Q}\textrm{ and }0<x<1\right\}$ is dense in $\textbf{R}$.\\
\\

Continuing we give a table of evaluations of certain $S(x)$, which as one can see depending from the discriminant $D$. The table can be verified using Theorem 13 below and using modular equations of the ''null'' Jacobi theta functions $\theta_3(q)$.\\ 
\\
Set
\begin{equation}
u:=S(x)\textrm{ and }v:=x^2,
\end{equation}
then\\
\textbf{1.} For $D=-4$
\begin{equation}
u=1
\end{equation}
\textbf{2.} For $D=-8$ 
\begin{equation}
-4u+4u^2+v=0
\end{equation}
\textbf{3.} For $D=-12$
\begin{equation}
-1-8u-18u^2+27u^4+16uv=0
\end{equation}
\textbf{4.} For $D=-16$
\begin{equation}
-64 u + 64 u^2 - 256 u^3 + 256 u^4 + 48 u v + 32 u^2 v + v^2=0
\end{equation}
\textbf{5.} For $D=-20$
\begin{equation}
-1 + 26 u - 275 u^2 + 1500 u^3 - 4375 u^4 + 6250 u^5 - 3125 u^6 - 
 256 u v + 256 u v^2=0
\end{equation}
\textbf{6.} For $D=-24$
$$
-6912 u^4 - 55296 u^5 - 124416 u^6 + 186624 u^8 - 1024 u v - 
 8192 u^2 v -
$$
$$
- 18432 u^3 v+ 
+ 6912 u^4 v + 69120 u^5 v + 62208 u^6 v + 
1280 u v^2 + 8480 u^2 v^2 +
$$
\begin{equation}
+16128 u^3 v^2 + 4320 u^4 v^2 - 
288 u v^3 + 112 u^2 v^3 + v^4=0
\end{equation}
\textbf{7.} For $D=-28$
$$
-1 - 48 u - 980 u^2 - 10976 u^3 - 72030 u^4 - 268912 u^5 - 
470596 u^6 + 823543 u^8 + 2144 u v - 
$$
\begin{equation}
-37632 u^2 v 
+21952 u^3 v + 
537824 u^5 v - 6144 u v^2 + 37632 u^2 v^2 + 4096 u v^3=0
\end{equation}
\textbf{8.} For $D=-32$
\begin{equation}
-16384 u + 16384 u^2 - 589824 u^3 + 589824 u^4 - 6291456 u^5 +6291456 u^6 -
$$
$$
-16777216 u^7 + 16777216 u^8 + 28672 u v - 
 221184 u^2 v + 933888 u^3 v - 327680 u^4 v +
$$
$$
 + 4456448 u^5 v + 
 1048576 u^6 v - 13568 u v^2 + 206336 u^2 v^2 - 372736 u^3 v^2 + 
 24576 u^4 v^2 +
$$
$$
+1216 u v^3 + 256 u^2 v^3 + v^4=0
\end{equation}
\textbf{9.} For $D=-36$
$$
-16384 u + 16384 u^2 - 589824 u^3 + 589824 u^4 - 6291456 u^5 + 
 6291456 u^6 -
$$
$$
 -16777216 u^7 + 16777216 u^8 + 28672 u v - 
 221184 u^2 v + 933888 u^3 v - 327680 u^4 v +
$$
$$
+4456448 u^5 v + 
 1048576 u^6 v - 13568 u v^2 + 206336 u^2 v^2 - 372736 u^3 v^2+
$$ 
\begin{equation}
+24576 u^4 v^2 + 1216 u v^3 + 256 u^2 v^3 + v^4=0
\end{equation}
\textbf{10.} For $D=-40$
$$
12800000 u^6 - 332800000 u^7 + 3520000000 u^8 - 19200000000 u^9 + 
 56000000000 u^{10}-
$$
$$ 
-80000000000 u^{11} + 40000000000 u^{12} - 
 262144 u v + 6815744 u^2 v - 72089600 u^3 v+
$$
$$
+ 393216000 u^4 v - 
 1153280000 u^5 v + 1785600000 u^6 v - 2163200000 u^7 v+
$$
$$  
+6080000000 u^8 v - 13600000000 u^9 v + 12000000000 u^{10} v + 
 589824 u v^2 -
$$
$$
-13631488 u^2 v^2 + 126156800 u^3 v^2 - 
588704000 u^4 v^2 + 1412480000 u^5 v^2 -
$$
$$
-1489600000 u^6 v^2 + 
 355200000 u^7 v^2 + 540000000 u^8 v^2 - 425984 u v^3+
$$
$$
+3899392 u^2 v^3 - 58668800 u^3 v^3 + 197984000 u^4 v^3 - 
293280000 u^5 v^3+
$$
$$
+10400000 u^6 v^3 + 102400 u v^4 + 
2920752 u^2 v^4 + 4463200 u^3 v^4 + 102000 u^4 v^4-
$$
\begin{equation}
-4200 u v^5+504 u^2 v^5 + v^6=0
\end{equation}
\\
\textbf{Theorem 12.}\\
If $(a,b,c)$ are in a $\theta-$class and  $D=-3,-7,-11,-19,-43,-67,-163$, then we have
\begin{equation}
\sum^{\infty}_{n,m=-\infty}q^{an^2+bnm+cm^2}=1+p\sum^{\infty}_{n=1}\left(\frac{n}{-D}\right)\frac{q^n}{1-q^n}\textrm{, }\forall |q|<1.
\end{equation}
where $p=6,2,2,2,2,2$ respectively.\\
\\
\textbf{Proof.}\\
Since for all $D$ mentioned in the state of the theorem we have $h(D)=1$, we get the validity of (65) in view of Lemma 1 and the notes of [9] section 4.3\\
\\
\textbf{Corollary 3.}\\ 
Let $r(n)$ be the number of representations of a positive integer $n$ into form
\begin{equation}
Q=ax^2+bxy+cy^2\textrm{, }x,y\in\textbf{Z}.
\end{equation}
When $D=-|D|=-3,-7,-11,-19,-43,-67,-163$, then 
\begin{equation}
r(n)=p\sum_{d|n}\left(\frac{d}{-D}\right)\textrm{, if }n\in\textbf{N},
\end{equation}
where $p=6,2,2,2,2,2,2$ respectively ($r(0)=1$).\\
\\
\textbf{Theorem 13.}\\
If $D=-|D|=b_1^2-4a_1c_1<0$, $a_1,c_1>0$ and exist $c_2$ positive integer such $b_1^2-4a_1c_1=-4a_2c_2$ with   $(a_1,b_1,c_1)\stackrel{\theta}{\equiv}(a_2,0,c_2)$, $a_2,c_2>0$, then 
\begin{equation}
\phi(a_1,b_1,c_1;q)=\theta_3\left(q^{a_2}\right)\theta_3\left(q^{c_2}\right)\textrm{, }|q|<1. 
\end{equation}
\\
\textbf{Proof.}\\
Given the triple $(a_1,b_1,c_1)$ if $D=-|D|=b_1^2-4a_1c_1=-4a_2c_2$ for some $a_2,c_2\geq 1$, then the triples $(a_1,b_1,c_1)$ and $(a_2,0,c_2)$ are equivalent and hence (68) holds.\\
\\ 
\textbf{Note 5.}\\The problem of the existence of such $\theta_3$ functions is solvable when $D\equiv 0(\textrm{mod} 4)$. There are no discriminants $D\equiv 2(\textrm{mod}4)$. From the cases $D=\pm 1(\textrm{mod}4)$ only $D\equiv 1(\textrm{mod}4)$ is occurring. Hence there remaining two cases: $D=8d+1$ and $D=8d+5$. The first of which occur iff $ac=2c_1$ and the second if $ac=2c_1-1$ and in both cases $c_1+d$ must be triangular number 'say' $t$. Hence given an odd discriminant, the values of $b$ and $ac$ are determined by a triangular number $t=\frac{k(k+1)}{2}$ as $b=2k+1$, $ac=\frac{(2k+1)^2-D}{4}$.\\ 
\\
\textbf{Examples 2.}\\
\textbf{i)} If $q=e^{-\pi\sqrt{r}}$, $r>0$, the theta function
\begin{equation}
\phi(q)=\sum^{\infty}_{n,m=-\infty}q^{n^2+2nm+3m^2}
\end{equation}
has discriminant $D=-8$, hence from (56)
\begin{equation}
-4 \left(\frac{\phi(q)}{\theta_3(q)^2}\right)^2+4\left(\frac{\phi(q)}{\theta_3(q)^2}\right)^4+k_r^2=0.
\end{equation}
Solving this equation and having in mind that $\theta_3(q)=\sum^{\infty}_{n=-\infty}q^{n^2}=\sqrt{\frac{2K}{\pi}}$, we get
\begin{equation}
\phi(q)=\frac{\theta_3^2(q)}{\sqrt{2}}\sqrt{1+k'_r}=\frac{\sqrt{2}}{\pi}K\sqrt{1+k'_r}
\end{equation}
\\
\textbf{ii)} If $q=e^{-\pi\sqrt{r}}$, $r>0$, the theta function
\begin{equation}
\phi(q)=\sum^{\infty}_{n,m=-\infty}q^{n^2+2nm+6m^2}
\end{equation}
has discriminant $D=-20$. Also  $0^2-4\cdot5=-20$ and $(1,2,6)\stackrel{\theta}{\equiv}(1,0,5)$. Hence from Theorem 13 we get
\begin{equation}
\phi(q)=\sum^{\infty}_{n,m=-\infty}q^{n^2+2nm+6m^2}=\theta_3\left(q\right)\theta_3\left(q^5\right)
\end{equation}
But $K(k_{25r})=m_5(r)K(k_r)=m_5K$. Hence
\begin{equation}
\sum^{\infty}_{n,m=-\infty}q^{n^2+2nm+6m^2}=\frac{2K\sqrt{m_5}}{\pi}
\end{equation}  
The equation for finding $m_5(r)=m_5$ is (see [12])
\begin{equation}
\left(5m_5-1\right)^5\left(1-m_5\right)=256\left(k_rk'_r\right)^2m_5
\end{equation}
By this way we can also get the case $D=-4$, $(a,b,c)\stackrel{\theta}{\equiv}(1,2,2)$ of Theorem 10.\\
\\
\textbf{iii)} The theta function
\begin{equation}
\phi(q)=\sum^{\infty}_{n,m=-\infty}q^{n^2+nm+2m^2}
\end{equation}
have discriminant $D=-7$ and if $\textbf{x}=(x,y)$, $A=\left(
\begin{array}{cc}
	2\textrm{ }1\\
	1\textrm{ }4
\end{array}
\right)$
clearly $A$ is an even symmetric $r\times r$, $r=2$ matrix with $\frac{1}{2}\textbf{x}^{t}A\textbf{x}=x^2+xy+2y^2$. The smallest positive integer $N$ such that $NA^{-1}$ is even matrix is $N=7$ (level). Hence $\theta(z)=\phi(e(z))$, $e(z)=e^{2\pi i z}$, $Im(z)>0$ is a modular form of weight $r/2=1$ such
\begin{equation}
\theta\left(\frac{az+b}{cz+d}\right)=\chi(d)(cz+d)\theta(z)\textrm{, for all }\left(
\begin{array}{cc}
	a\textrm{ }b\\
	c\textrm{ }d
\end{array}
\right)\in \Gamma_0(7),
\end{equation} 
$\chi(d)=\left(\frac{D_1}{d}\right)=\left(\frac{-7}{d}\right)$ (Kronecker symbol) since $D_1=(-1)^{r/2} \textrm{det}(A)=-7$. Hence $\theta(z)^2$ is a modular form of weight 2 in $\Gamma_0(7)$ and $\textrm{dim}E_2\left(\Gamma_0(7)\right)=1$, $\textrm{dim}S_2\left(\Gamma_0(7)\right)=0$. Hence we can expand $\theta(z)^2$ into Eisenstein series. We find
\begin{equation}
\left(\sum^{\infty}_{n,m=-\infty}q^{n^2+nm+2m^2}\right)^2=\frac{-1}{6}\left(E_2(z)-7E_2(7z)\right)\textrm{, }q=e(z)\textrm{, }Im(z)>0,
\end{equation}
where
\begin{equation}
E_{k}(z)=1-\frac{2k}{B_{k}}\sum^{\infty}_{n=1}\sigma_{k-1}(n)q^n
\end{equation}
is the $k$-weight Eisenstein series and  $\sigma_{\nu}(n):=\sum_{d|n}d^{\nu}$.\\
Moreover we define the arithmetical function 
\begin{equation}
\sigma^{*}(n):=\sum_{\scriptsize
\begin{array}{cc}
	d|n\\
	(d,7)=1\normalsize
\end{array}
}d.
\end{equation}
Clearly $\sigma^{*}(n)$ is a multiplicative function and for every prime $p$ satisfies the identity
\begin{equation}
\sigma^{*}\left(p^{n+1}\right)=\sigma^{*}\left(p\right)\sigma^{*}\left(p^{n}\right)-p^{2-1}\cdot\sigma^{*}\left(p^{n-1}\right)
\end{equation}
Hence the series
\begin{equation}
1+\sum^{\infty}_{n=1}\sigma^{*}\left(n\right)q^n
\end{equation} 
is a Hecke form of weight 2. Observe that
\begin{equation}
\frac{T_n\theta^2(z)}{\theta^2(z)}=\sigma^{*}\left(n\right)
\end{equation}
and consequently
\begin{equation}
\left(\sum^{\infty}_{n,m=-\infty}q^{n^2+nm+2m^2}\right)^2=1+\sum^{\infty}_{n=1}\sigma^{*}(n)q^n.
\end{equation}
By this way, we get the cases $D=-3,-7$ of Theorem 10 and Corollary 2 for $(a,b,c)\stackrel{\theta}{\equiv}(1,1,1)$, $(a,b,c)\stackrel{\theta}{\equiv}(1,1,2)$.\\
\\
\textbf{iv)} As in (iii) we have $D=-11$,
\begin{equation}
\phi(q)=\sum^{\infty}_{n,m=-\infty}q^{n^2+nm+3m^2}
\end{equation}
and $A=\left(
\begin{array}{cc}
	2\textrm{ }1\\
	1\textrm{ }6
\end{array}
\right)$, $A^{-1}=\left(
\begin{array}{cc}
	\frac{6}{11}\textrm{ }\frac{-1}{11}\\
	\frac{-1}{11}\textrm{ }\frac{2}{11}
\end{array}
\right)$. Hence the level of $\phi(q)$ is $N=11$, $r=2$. Therefore 
\begin{equation}
\theta(z)=\phi(q)\textrm{, }q=e(z)\textrm{, }Im(z)>0
\end{equation}
is a modular form of weight 1 on $M_1\left(\Gamma_0(11),\chi\right)$. As before $D_1=(-1)^{r/2}\textrm{det}(A)=-11$,
\begin{equation}
\theta\left(\frac{az+b}{cz+d}\right)=\chi(d)(cz+d)\theta(z)
\end{equation}
and $\chi(d)=\left(\frac{-11}{d}\right)$, for all $\left(
\begin{array}{cc}
	a\textrm{ }b\\
	c\textrm{ }d
\end{array}
\right)\in \Gamma_0(11)$. Therefore we obtain that $\theta(z)^2$ is a modular form of wight 2 in $\Gamma_0(11)$; $\textrm{dim}E_2(\Gamma_0(11))=\textrm{dim}S_2(\Gamma_0(11))=1$. A base of $E_2(\Gamma_0\left(11)\right)$ is $E_2(z)-11E_2(11z)$ and a base for $S_2\left(\Gamma_0(11)\right)$ is $\eta(z)^2\eta(11z)^2$. Hence we easily conclude that
\begin{equation}
\left(\sum^{\infty}_{n,m=-\infty}q^{n^2+nm+3m^2}\right)^2=-\frac{1}{10}\left(E_2(z)-11E_2(11z)\right)+\frac{8}{5}\eta(z)^2\eta(11z)^2.
\end{equation} 
\\
\textbf{Note 6.}\\
1) It is clear that if $a,b,c\in\textbf{N}$ with $\textrm{GCD}(a,b,c)=1$ and $D=-|D|$ is odd, then the level $N$ of the modular form 
\begin{equation}
\sum^{\infty}_{n,m=-\infty}q^{an^2+bnm+cm^2}\textrm{, }q=e(z),
\end{equation}
is $N=-D$.\\
2) If $l\geq 2$ and $l|N$, then $E_2(z)-lE_2(lz)\in M_2(\Gamma_0(N))$\\  
3) $M_k(\Gamma_0(l))\subseteq M_k(\Gamma_0(m))$ if $l|m$.\\
\\
\textbf{Lemma 2.}\\
Let $l$ be a positive integer and consider the eta-product $f(z)$ defined as
\begin{equation}
f(z)=\prod_{\delta |l}\left(\eta(\delta z)\right)^{r_{\delta}}.
\end{equation} 
Let also
\begin{equation}
k=\frac{1}{2}\sum_{\delta |l}r_{\delta}\textrm{ and }s=\prod_{\delta|l}\delta^{|r_{\delta}|}.
\end{equation}
Then if $f(z)$ satisfies all conditions $1,2,3,4,5$ below\\
1) $k$ in an even integer\\
2) $s$ is the square of an integer\\
3) 
$$
\sum_{\delta|l}\delta r_{\delta}\equiv 0(\textrm{mod}24)
$$
4) 
$$
\sum_{\delta|l}\frac{l}{\delta}r_{\delta}\equiv 0(\textrm{mod}24)
$$
5) 
$$
\sum_{\delta|l}\textrm{GCD}(d,\delta)^2\frac{r_{\delta}}{\delta}\geq 0\textrm{, for all }d|l,
$$ 
$f$ is a modular form in $M_k(\Gamma_0(l))$.\\
\\
\textbf{Proof.}\\
See article [8] Lemma 2.1 and the related references.\\
\\   
\textbf{v)} If $(a,b,c)\stackrel{\theta}{\equiv}(1,1,4)$, $D=-15$, $N=15$, then $\textrm{dim}E_2(\Gamma_0(15))=3$, $\textrm{dim}S_2(\Gamma_0(15))=1$. A base for the cusp forms is $\eta(z)\eta(3z)\eta(5z)\eta(15z)$ and a base for the Eisenstein subspace consist of $E_2(z)-3E_2(3z)$, $E_2(z)-5E_2(5z)$, $E_2(z)-15E_2(15z)$. We can easily verify that
$$
\left(\sum^{\infty}_{n,m=-\infty}q^{n^2+nm+4m^2}\right)^2=\frac{1}{12}\left(E_2(z)-3E_2(3z)\right)-\frac{1}{12}\left(E_2(z)-5E_2(5z)\right)+
$$
\begin{equation}
+\frac{1}{12}\left(E_2(z)-15E_2(15z)\right)-2\eta(z)\eta(3z)\eta(5z)\eta(15z).
\end{equation}
\\
\textbf{vi)} For the case $(a,b,c)\stackrel{\theta}{\equiv}(1,1,7)$, $-D=27=N$, we have $\textrm{dim}E_2(\Gamma_0(27))=5$, $\textrm{dim}S_2(\Gamma_0(27))=1$ and (see [8]):  
$$
\left(\sum^{\infty}_{n,m=-\infty}q^{n^2+nm+7m^2}\right)^2=C_1(E_2(z)-3E_2(3z))+C_2(E_2(z)-9E_2(9z))+
$$
\begin{equation}
+C_3(E_2(z)-27E_2(27z))+C_4E_{2,\chi}(z)+C_5E_{2,\chi}(3z)+C_6\Delta(z),
\end{equation}
where
\begin{equation}
\Delta(z)=\eta(3z)^2\eta(9z)^2
\end{equation} 
and 
\begin{equation}
E_{k,\chi}(z)=\sum^{\infty}_{n=1}\chi(n)\sigma_{k-1}(n)q^n\textrm{, }q=e(z)\textrm{, }Im(z)>0.
\end{equation}
In (93) the constants are $C_1=-\frac{2}{27}$, $C_2=\frac{2}{27}$, $C_3=-\frac{1}{18}$, $C_4=C_5=0$, $C_6=\frac{8}{3}$.\\ 
\\
\textbf{vii)} For the case $-D=35=N$, $(a,b,c)\stackrel{\theta}{\equiv}(1,1,9)$, we have  $\textrm{dim}S_2(\Gamma_0(35))=3$, $\textrm{dim}E_2(\Gamma_0(35))=3$.\\Hence using Lemma 2 we find the cusp forms 
$$
\frac{\eta(z)^3\eta(35z)^3}{\eta(5z)\eta(7z)}\textrm{, }\eta(z)\eta(5z)\eta(7z)\eta(35z),
$$
\begin{equation}
\eta(5z)^2\eta(7z)^2\textrm{, }\eta(z)^2\eta(35z)^2.
\end{equation}
The form
\begin{equation}
\frac{\eta(5z)^3\eta(7z)^3}{\eta(z)\eta(35z)}
\end{equation}
is not a cusp form, since it not vanish at $i\infty$.  
Hence assuming all the possible pairs we find
$$
\left(\sum^{\infty}_{n,m=-\infty}q^{n^2+nm+9m^2}\right)^2=-
\frac{1}{2}(E_2(z)-35E_2(35z))
-16\frac{\eta(5z)^3\eta(7z)^3}{\eta(z)\eta(35z)}+
$$
\begin{equation}
+8\eta(5z)^2\eta(7z)^2\textrm{, }q=e(z)\textrm{, }Im(z)>0.
\end{equation} 
 
Also hold relations like
\begin{equation}
\phi(1,1,3;q)=1+2\sum^{11}_{j=1}\left(\frac{j}{11}\right)\sum^{\infty}_{n=0}\frac{q^{11n+j}}{1-q^{11n+j}}
\end{equation}
and
\begin{equation}
\phi(1,1,7;q)=1+2\sum^{27}_{j=1}A_j\sum^{\infty}_{n=0}\frac{q^{27n+j}}{1-q^{27n+j}},
\end{equation}
where
\begin{equation}
A_l=\left\{
\begin{array}{cc}
  u_l\textrm{, if }3|l\\
  \left(\frac{l}{27}\right)\textrm{, else }
\end{array}\right\},
\end{equation}
with  $u_k=0,-1,1,3,-1,1,-3-1,1$, if $k\equiv 0,3,6,9,12,15,18,21,24\left(\textrm{mod}27\right)$, resp.
Hence
\begin{equation}
\sum^{\infty}_{n,m=-\infty}q^{n^2+nm+3m^2}=1+2\sum^{\infty}_{n=1}\left(\sum_{d|n}\left(\frac{d}{11}\right)\right)q^n
\end{equation}
and
\begin{equation}
\sum^{\infty}_{n,m=-\infty}q^{n^2+nm+7m^2}=1+2\sum^{\infty}_{n=1}\left(\sum_{d|n}A_d\right)q^n.
\end{equation}

\section{The modularity of $\sum^{\infty}_{n,m=-\infty}q^{an^2+bnm+cm^2}$}

Suppose $q=e^{2\pi i z}=e(z)$, with $Im(z)>0$. Also let $a>0$, $b^2-4ac<0$ and 
\begin{equation}
f_0(x,y)=\exp\left[-\frac{\sqrt{t}}{4\pi}(ax^2+bxy+cy^2)\right].
\end{equation}
If we define the Fourier transform of $f(x,y)$ as
\begin{equation}
\widehat{f}(s,w)=\frac{1}{2\pi}\int\int_{\textbf{\small R\normalsize}^2}f(x,y)e^{ixs+iyw}dxdy
\end{equation}
then
\begin{equation}
\widehat{f_0}(s,w)=\frac{8\pi^2}{\sqrt{-D}\sqrt{t}}\exp\left[\frac{4\pi}{D\sqrt{t}}\left(as^2-bsw+cw^2\right)\right].
\end{equation} 
The Poisson summation formula read as
\begin{equation}
\sum^{\infty}_{n,m=-\infty}f(2\pi n,2\pi m)=\frac{1}{4\pi^2}\sum^{\infty}_{n,m=-\infty}\widehat{f}(n,m).
\end{equation} 
Hence if we denote
\begin{equation}
\phi(q):=\sum^{\infty}_{n,m=-\infty}q^{an^2+bnm+cm^2}\textrm{, }q=e(z)
\end{equation}
and
\begin{equation}
\psi(z):=\phi\left(e\left(\frac{z}{\sqrt{-D}}\right)\right)=\phi\left(q^{1/\sqrt{-D}}\right),
\end{equation}
then in view of the Poisson summation formula we have
$$
\sum_{(n,m)\in\textbf{\small Z\normalsize}^2}\exp\left(-\sqrt{t}\pi\left(an^2+b nm+cm^2\right)\right)=
$$
\begin{equation}
=\frac{2}{\sqrt{-D}\sqrt{t}}\sum_{(n,m)\in\textbf{\small Z\normalsize}^2}\exp\left(\frac{2}{iD\sqrt{t}}2\pi i\left(an^2+bnm+cm^2\right)\right)
\end{equation}
and get the next\\
\\
\textbf{Theorem 14.}\\
If $a>0$ and $D=b^2-4ac<0$, the function $\psi(z)$ of (109) is homomorphic in the upper half plane, periodic with period $\sqrt{-D}$ and 
\begin{equation}
\psi\left(\frac{-1}{z}\right)=\frac{z}{i}\psi(z).
\end{equation}
Also $\phi(e(z))$ of (108) is $1-$periodic and
\begin{equation}
\phi\left(e\left(\frac{-1}{|D|z}\right)\right)=\frac{z \sqrt{|D|}}{i}\phi\left(e(z)\right)
\end{equation}
\\
\textbf{Note 7.}\\
\textbf{1)} Assume the Dedekind eta function
\begin{equation}
\eta(z)=q^{1/24}\prod^{\infty}_{n=1}\left(1-q^n\right)\textrm{, }q=e(z)\textrm{, }Im(z)>0. 
\end{equation}
When $Im(z)>0$ the functions
\begin{equation}
\upsilon(z)=\eta(\alpha z)\eta(\beta z)\textrm{, with }\alpha+\beta=24\textrm{, }\alpha,\beta\in\textbf{N}^{*}
\end{equation} 
satisfy the following functional equation
\begin{equation}
\upsilon\left(\frac{-1}{\delta_1 z}\right)=\frac{z\sqrt{\delta_1}}{i}\upsilon(z)\textrm{, }\delta_1=\alpha \beta
\end{equation}
which is the same to that of $\phi(e(z))$ (relation (112)). The functions $\upsilon(z)$ are modular forms (cusp forms) in the group $M_1(\Gamma_0(\alpha\beta))$.\\
\textbf{2)} The functions $\phi(e(z))=\phi(a,b,c;e(z))$ satisfy (112) also are $1-$periodic and belong to $M_1(\Gamma_0(-D),\chi)$. More precisely when $-D$ is odd positive integer they are modular forms of weight 1 and of character $\chi$ and satisfy the functional equation
\begin{equation}
\phi\left(e\left(\frac{az+b}{cz+d}\right)\right)=\chi(d)(cz+d)\phi\left(e(z)\right),
\end{equation}
where $\chi(d)=\left(\frac{D}{d}\right)$, for all $\left(
\begin{array}{cc}
	a\textrm{ }b\\
	c\textrm{ }d
\end{array}
\right)\in \Gamma_0(-D)$.\\
If 
\begin{equation}
f(z)=\phi(e(z))=\sum^{\infty}_{n=0}a_f(n)e(nz),
\end{equation}
then $a_f(n)$, is the number of representations of the positive integer $n$ in the form $ax^2+bxy+cy^2$. If $G(s):=(2\pi)^{-s}\Gamma(s)$, then 
$$
\Lambda_f(s)=\int^{+\infty}_{0}\phi(e(it))t^{s-1}dt=G(s)\sum^{\infty}_{n=1}\frac{a_f(n)}{n^s}=
$$
\begin{equation}
=G(s)\sum_{(n,m)\in\textbf{\scriptsize{Z}}^2-\{(0,0)\}}\frac{1}{(an^2+bnm+cm^2)^s}
\end{equation}  
and
\begin{equation}
\Lambda_f(s)=|D|^{1/2-s}\Lambda_f(1-s).
\end{equation}
Equation (119) is a necessary condition for $f$ to belong to $M_1(\Gamma_0(-D),\chi)$.\\
\\
\textbf{Theorem 15.}\\
Let $f(z)=\phi(a,b,c;q)$ and $a(n)$ be such that
$$
f(z)^2=1+\sum^{\infty}_{n=1}a(n)q^n\textrm{, }q=e(z)\textrm{, }Im(z)>0,
$$
then
$$
\Lambda(s)=\int^{+\infty}_{0}\phi(e(it))^2t^{s-1}dt=G(s)\sum^{\infty}_{n=1}\frac{a(n)}{n^s}
$$
and
\begin{equation}
\Lambda(s)=|D|^{1-s}\Lambda(2-s)
\end{equation}
and $f(z)^2$ belongs to $M_2(\Gamma_0\left(-D)\right)$.\\
\\
\textbf{Proof.}\\
We have from relation (112):
\begin{equation}
\phi\left(e\left(\frac{i}{|D|t}\right)\right)^2=t^2|D|\phi\left(e(it)\right)^2.
\end{equation}
Hence
$$
\Lambda(s)=\int^{+\infty}_{0}\phi(e(it))^2t^{s-1}dt=|D|\int^{+\infty}_{0}\phi\left(e\left(\frac{i}{|D|t}\right)\right)^2\frac{1}{t^2|D|^2}t^{s-1}dt.
$$
Making the change of variable $w=\frac{1}{|D|t}$, we have
$$
\Lambda(s)=|D|\int^{0}_{+\infty}\phi\left(e\left(iw\right)\right)^2w^2\left(\frac{1}{|D|w}\right)^{s-1}\frac{-dw}{|D|w^2}=
$$
$$
=|D|^{1-s}\int^{+\infty}_{0}\phi\left(e\left(iw\right)\right)^2w^{1-s}dw=|D|^{1-s}\Lambda(2-s).
$$
\\
\textbf{Theorem 16.}\\
The function $f(z)=\phi(a,b,c,e(z))^{2\nu}$, with $a,c,\nu\in\textbf{N}$ and $b$ integer such that $D=b^2-4ac$ is negative, is a modular form of weight $2\nu$ and belongs to $M_{2\nu}(\Gamma_0(-D))$. It holds $M_{2\nu}(\Gamma_0(-D))=E_{2\nu}(\Gamma_0(-D))\oplus S_{2\nu}(\Gamma_0(-D))$. The space $E_{2\nu}(\Gamma_0(-D))$ have finite base constructed from linear combination of Eisenstein series of weight $2\nu$ in $\Gamma_0(-D)$. The second space $S_{2\nu}(\Gamma_0(-D))$ is also finite dimensional and is linear combination of cusp forms of $\Gamma_0(-D)$ and weight $2\nu$.\\ 
\\
\textbf{Theorem 17.} (Hecke-Weil)\\
If the function $f(z)$ is homomorphic in the upper half plane $Im(z)>0$ with Fourier expansion 
\begin{equation}
f(z)=\sum^{\infty}_{n=0}a_f(n)e(nz)\textrm{, in } Im(z)>0,
\end{equation}
and
\begin{equation}
f(z)=O(Im(z)^{-\nu})\textrm{, }Im(z)\rightarrow 0\textrm{, }\nu>0,
\end{equation}
then for positive integers $|D|$, $k/2$ the following are equivalent:\\
\textbf{1.} 
\begin{equation}
f\left(\frac{-1}{|D|z}\right)=\left(\frac{\sqrt{|D|}z}{i}\right)^kf(z)\textrm{, }\forall z: Im(z)>0
\end{equation}
\textbf{2.}
The coefficients $a_f(n)$ form a function
\begin{equation}
\Lambda^{*}_f(s):=\left(\frac{\sqrt{|D|}}{i}\right)^k\int^{+\infty}_{0}f(it)t^{s-1}dt=G(s)\left(\frac{\sqrt{|D|}}{i}\right)^k\sum^{\infty}_{n=1}\frac{a_f(n)}{n^s},
\end{equation}
(where $G(s)=(2\pi)^{-s}\Gamma(s)$) with the properties\\
\textbf{a)}
\begin{equation}
\Lambda^{*}_f(s)=i^{k}\Lambda^{*}_f(k-s)
\end{equation}
and\\
\textbf{b)}
\begin{equation}
\Lambda^{*}_f(s)+\frac{a_f(0)}{s}+\frac{a_f(0)i^k}{k-s}
\end{equation}
is homomorphic on the $s-$plane and bounded on any vertical strip.\\
\textbf{3.} $f(z)$ is a modular form of weight $k$ in $M_k\left(\Gamma_0(|D|)\right)$.\\
\\
\textbf{Example 3.}\\
\textbf{i)} Consider the function
\begin{equation}
f(z):=\phi(1,1,1;q)^2=\left(\sum^{\infty}_{n,m=-\infty}q^{n^2+nm+m^2}\right)^2\textrm{, }q=e(z)\textrm{, }Im(z)>0.
\end{equation}
Clearly $f(z)$ belongs to $M_2(\Gamma_0(3))=E_2(\Gamma_0(3))\oplus S_2(\Gamma_0(3))$ and $\textrm{dim}E_2(\Gamma_0(3))=1$, $\textrm{dim}S_2(\Gamma_0(3))=0$. An element of $E_2(\Gamma_0(3))$ is $E_2(z)-3E_2(3z)$. Hence easily
\begin{equation}
\left(\sum^{\infty}_{n,m=-\infty}q^{n^2+nm+m^2}\right)^2=-\frac{1}{2}\left(E_2(z)-3E_2(3z)\right)
\end{equation}
\\ 
\textbf{ii)} The function $f(z)=\phi(1,1,1,q)^4$ is a modular form of weight 4 and belongs to $M_4(\Gamma_0(3))$. This space can be spliten to $E_{4}(\Gamma_0(3))$ and $S_{4}(\Gamma_0(3))$. The later space of cusps forms has dimension zero while $\textrm{dim}E_4(\Gamma_0(3))=2$. It must hold that 
\begin{equation}
\phi(1,1,1;q)^4=C_1E_4(z)+C_2\left(E_{2}(z)-E_2\left(3z\right)\right)^2
\end{equation}  
since $(E_2(z)-3E_2(3z))^2$ and $E_4(z)$ are linear independent and belong to $E_4(\Gamma_0(3))$. Hence comparing coefficients we find
\begin{equation}
\left(\sum^{\infty}_{n,m=-\infty}q^{n^2+nm+m^2}\right)^4=\frac{1}{10}E_4(z)+\frac{9}{10} E_4(3z).
\end{equation}
Also if we consider the form
\begin{equation}
T=\sum^{4}_{i=1}(x_i^2+x_iy_i+y_i^2),
\end{equation}
then the number of representations of a given $n$ to $T$ is
\begin{equation}
r_T(n)=24\sigma_3(n)+216\sigma_{3}\left(\frac{n}{3}\right)\textrm{, }r_T(0)=1
\end{equation}
The sigma function is zero when 3 not divides $n$.\\

In general hold many relations easy to check, especially those having no cusp forms like
\begin{equation}
\left(\sum^{\infty}_{n,m=-\infty}q^{n^2+nm+m^2}\right)^2=-\frac{1}{2}E_{2}(z)+\frac{3}{2}E_2(3z)
\end{equation}
and
\begin{equation}
\left(\sum^{\infty}_{n,m=-\infty}q^{n^2+nm+2m^2}\right)^2=-\frac{1}{6}E_2(z)+\frac{7}{6}E_2(7z)
\end{equation}
\\

We give now a parametric family of functions, similar to Eisenstein, that are also modular forms.\\
\\
\textbf{Theorem 18.}\\
For $k=2,4,6\ldots$ the functions
\begin{equation}
H_{k}(z)=\sum_{(n,m)\in\textbf{A}}\frac{(-1)^n}{(n+mz)^{k}},
\end{equation}
where $\textbf{A}=\textbf{Z}\times\textbf{Z}-\{(0,0)\}$, are modular forms of weight $k$ and in particular belong to  $M_k(\Gamma_1(2))$. Further if $q=e_1(z):=e^{i\pi z}$, $Im(z)>0$, then
\begin{equation}
H_{k}(z)=\frac{\pi^{k}}{k!}\left[\left(2-2^k\right)|B_{k}|+4ki^kF_{k}(q)\right],
\end{equation}
where
\begin{equation}
F_{k}(q):=\sum^{\infty}_{n=1}\sigma^{*}_{k-1}(n)q^{n}
\end{equation}
and $\sigma^{*}_{\nu}(n):=\sum_{d-odd,d|n}d^{\nu}$.\\
\\
\textbf{Proof.}\\
i) The function is clearly a modular form since the action of transformation $a\equiv1(2)$, $b\equiv0(2)$, $c\equiv 0(2)$, $d\equiv 1(2)$, with $ad-bc=1$ gives immediately using (136)
$$
H_{k}\left(\frac{az+b}{cz+d}\right)=(cz+d)^{k}H_{k}(z).
$$\\
ii) Set $q=e^{i\pi z}$ with $Im(z)>0$. Using the formulas 
\begin{equation}
\pi \csc(\pi z)=\sum^{\infty}_{n=-\infty}\frac{(-1)^n}{z+n}
\end{equation}
and 
$$
\csc(\pi z)=\frac{2i}{q-q^{-1}},
$$ 
we obtain differentiating (139) $k$ times
\begin{equation}
\sum^{\infty}_{n=-\infty}\frac{(-1)^n}{(z+n)^k}=2\frac{(-\pi i)^k}{\Gamma(k)}\sum^{\infty}_{n=1}(2n-1)^{k-1}q^{2n-1}
\end{equation}
Set where $z\rightarrow m\tau$ and suppose that $k$ is even. Sum with respect to $m$ to get
\begin{equation}
\sum^{\infty}_{m=1}\sum^{\infty}_{n=-\infty}\frac{(-1)^n}{(n+mz)^k}=2\frac{(\pi i)^k}{\Gamma(k)}\sum^{\infty}_{m,n=1}(2n-1)^{k-1}q^{(2n-1)m}
\end{equation}
Lastly we can write for $k=1,2,\ldots$ the identity
\begin{equation}
\sum_{(n,m)\in\textbf{A}}\frac{(-1)^n}{(n+mz)^{2k}}=2\sum^{\infty}_{n=1}\frac{(-1)^n}{n^{2k}}+4\frac{(\pi i)^{2k}}{\Gamma(2k)}F_{2k}(q)
\end{equation} 
where 
$$
F_{2k}(q)=\sum^{\infty}_{n=1}\sigma^{*}_{2k-1}(n)q^n
$$
The result follows from (see [14] ch.23)
$$
\sum^{\infty}_{n=1}\frac{(-1)^n}{n^s}=(2^{1-s}-1)\zeta(s)\textrm{, }s\geq 1
$$
and
$$
\zeta(2s)=\frac{(2\pi)^{2s}}{2(2s)!}|B_{2s}|\textrm{, }s=1,2,\ldots
$$
$B_s$ are the Bernoulli numbers.\\

It is known that if $q=e(z)$, then $\theta(z)=\sum^{\infty}_{n=-\infty}q^{n^2/2}$ belongs to the modular group $M_{\frac{1}{2}}(\Gamma_1(4))$. Hence if $a,d\equiv 1(4)$, $b,c\equiv 0(4)$, $ad-bc=1$, then $\theta$ is a modular form of weight $1/2$ and satisfy the following functional equation:
\begin{equation}
\theta\left(\frac{az+b}{cz+d}\right)=(cz+d)^{1/2}\theta(z).
\end{equation} 
Also from the above theorem we have that $H_2(z)$ is element of $M_2\left(\Gamma_0(2)\right)$, but it is $M_2\left(\Gamma_0(2)\right)\subset M_2\left(\Gamma_0(4)\right)$. Hence the set $M_2\left(\Gamma_0(4)\right)$ have elements such  $H_2(z)$, $H^{*}_2(z)$, $E_2(z)-2E_2(2z)$, $E_2(z)-4E_2(4z)$, $\theta(z)^4$, where the sign $*$ means $H^{*}_{k}(z)=H_k(z+1)$. Hence 
\begin{equation}
\theta(z)^4+C_1H_{2}(z)+C_2H^{*}_2(z)=0, 
\end{equation}
where $C_1=\frac{4}{\pi^2}$, $C_2=\frac{2}{\pi^2}$.\\

Another example is the power $\theta(z)^8$ which is a modular form of weight 4 in $\Gamma_1(4)$ and can be constructed using the forms $H_4$,$H^{*}_4$ and $E_4$ (here $E_{2k}(z)$ are the known  Eisenstein series of weight $2k$). More detailed we have
\begin{equation}
\theta(z)^8+\frac{24}{\pi^4}H_4(z)+\frac{24}{\pi^4}H_4^{*}(z)-\frac{1}{15}E_4(z)=0
\end{equation}
As application of the above formula holds the next\\
\\
\textbf{Theorem 19.}\\
The number of representations of a number $n\in\textbf{N}$ as the sum of 8 squares is
\begin{equation}
r_8(n)=16\sigma_{3}(n)-32\cdot\epsilon_n\sum_{d|n,d-odd}d^3
\end{equation}
where $\epsilon_n=1$ if $n$ is even and 0 else.\\
\\
\textbf{Corollary 4.}\\
Every odd number can be written as a sum of 8 squares.\\
\\ 
\textbf{Theorem 20.}\\
We can square (145) to get $\theta(z)^{16}$ and hence $r_{16}(n)$:
$$
r_{16}(n)=32\sigma_3(n)-64\epsilon_n\sigma^{*}_3(n)+256\sum^{n}_{m=1}\sigma_3(m)\sigma_3(n-m)-
$$
$$
-512\sum^{n}_{m=1}\left(1+(-1)^{n-m}\right)\sigma_3(m)\sigma_3^{*}(n-m)+
$$
\begin{equation}
+512\sum^{n}_{m=1}\left(\epsilon_n+(-1)^{n-m}\right)\sigma^{*}_3(m)\sigma^{*}_3(n-m),
\end{equation}
where $\sigma_3(0)=0$ and $\sigma^{*}_3(0)=0$.

\section{The Ramanujan expansions}

For $|q|<1$, the Rogers Ramanujan continued fraction (RRCF) is defined as
\begin{equation}
R(q):=\frac{q^{1/5}}{1+}\frac{q}{1+}\frac{q^2}{1+}\frac{q^3}{1+}\ldots
\end{equation}
If $q=e^{-\pi\sqrt{r}}$, $r>0$ then
\begin{equation}
f(-q):=\prod^{\infty}_{n=1}\left(1-q^n\right)
\end{equation}
is the Ramanujan eta function and
\begin{equation}
\eta(q):=q^{1/24}\prod^{\infty}_{n=1}\left(1-q^n\right)\textrm{, when }|q|<1,
\end{equation}
is the Dedekind eta function. We also set
\begin{equation}
\eta_1(z):=q^{1/24}\prod^{\infty}_{n=1}\left(1-q^n\right)\textrm{, }q=e(z):=e^{2\pi i z}\textrm{, }Im(z)>0
\end{equation}
and
\begin{equation}
A_r:=\frac{f\left(-q^2\right)^6}{q^2f\left(-q^{10}\right)^6}=\frac{\eta_1\left(\sqrt{-r}\right)^6}{\eta_1\left(5\sqrt{-r}\right)^6}=R\left(q^2\right)^{-5}-11-R\left(q^2\right)^5,
\end{equation}
for $q=e^{-\pi\sqrt{r}}$, $r>0$. The function
\begin{equation}
\Pi(r):=3\sqrt[3]{2k_{4r}}\cdot { }_2F_1\left(\frac{1}{3},\frac{1}{6};\frac{7}{6};k_{4r}^2\right)
\end{equation}
is called Carty's function is related to famous Carty's problem, and satisfies the functional equation
\begin{equation}
\Pi(r)+\Pi\left(\frac{1}{r}\right)=C_0\textrm{, }\forall r>0, 
\end{equation}
with $C_0=2^{-4/3}\pi^{-1}\Gamma(1/3)^{3}\sqrt{3}$. Its first derivative is:
\begin{equation}
\frac{d}{dr}\Pi(r)=-\pi\frac{\eta_1\left(i\sqrt{r}\right)^4}{\sqrt{r}}.
\end{equation}
The function $k_r$ is the elliptic singular modulus i.e. the solution of the equation 
\begin{equation}
\frac{K\left(\sqrt{1-k_r^2}\right)}{K\left(k_r\right)}=\sqrt{r}\textrm{, }r>0
\end{equation}
where
\begin{equation} K(x)=\frac{\pi}{2}\cdot{}_2F_1\left(\frac{1}{2},\frac{1}{2};1;x^2\right)\textrm{, }|x|<1
\end{equation}
is the complete elliptic integral of the first kind. 
It is known that $k_r$ is algebraic when $r$ is positive rational.\\
The following are the traditional definitions of the complete elliptic integrals of the first and second kind respectively
\begin{equation}
K(x)=\int^{\pi/2}_{0}\frac{dt}{\sqrt{1-x^2\sin^2(t)}} \textrm{ and } E(x)=\int^{\pi/2}_{0}\sqrt{1-x^2\sin^2(t)}dt
\end{equation}
The elliptic alpha function is (see [12]): 
\begin{equation}
\alpha(r):=\frac{\pi}{4K^2(k_r)}-\sqrt{r}\left(\frac{E(k_r)}{K(k_r)}-1\right).
\end{equation}
The function $\alpha(r)$ has also the property of being algebraic when $r$ is positive rational.\\

We prove the next\\
\\
\textbf{Theorem 21.}\\
If $r>0$, then 
\begin{equation}
\int^{+\infty}_{A_r}\frac{dt}{t^{1/6}\sqrt{125+22t+t^2}}=\Pi(r)=\frac{1}{\sqrt[3]{4}}B\left(k_{4r}^2,1/6,2/3\right),
\end{equation}
where $B(x,a,b)=\int^{x}_{0}t^{a-1}(1-t)^{b-1}dt$ is the incomplete beta function.\\
\\
\textbf{Proof.}\\
If $q=e^{-\pi\sqrt{r}}$, $r>0$, we define 
\begin{equation}
T_{r}:=1+6\sum^{\infty}_{n=1}\frac{nq^{2n}}{1-q^{2n}}-30\sum^{\infty}_{n=1}\frac{nq^{10n}}{1-q^{10n}}=\frac{-1}{4}\left(P\left(q^2\right)-5P\left(q^{10}\right)\right), 
\end{equation}
with $P(q)=1-24\sum^{\infty}_{n=1}\frac{nq^n}{1-q^n}$ being the familiar Eisenstein series, one can easily see that
\begin{equation}
T_{r}=\frac{\sqrt{r}}{\pi A_r}\cdot\frac{dA_r}{dr}.
\end{equation}
Ramanujan have found and later proven by Berndt (see [5] pg.464 Entry 4(ii)) that
\begin{equation}
1+6\sum^{\infty}_{n=1}\frac{nq^{2n}}{1-q^{2n}}-30\sum^{\infty}_{n=1}\frac{nq^{10n}}{1-q^{10n}}=\frac{\left(x^{2}+22 q^2 x y +125 q^4y^2\right)^{1/2}}{\left(xy\right)^{1/6}}, 
\end{equation}
where $x=f\left(-q^2\right)^6$ and $y=f\left(-q^{10}\right)^6$. But also it is
\begin{equation}
y=\frac{x}{q^2 A_r}.
\end{equation}
Setting this $y$ in relation (163) we get
$$
T_r=q^{1/3} x^{2/3}A^{-5/6}_r\sqrt{125+22A_r+A^{2}_r}.
$$
Hence from this last equation and relations (162),(152),(149),(150), we get
\begin{equation}
\frac{dA_r}{dr}=\pi\frac{\eta\left(i\sqrt{r}\right)^4}{\sqrt{r}}A_r^{1/6}\sqrt{125+22 A_{r}+A_{r}^2}.
\end{equation}
Integrating equation (165) we get the result.\\

Ramanujan have evaluated all 
\begin{equation}
T_{p,r}=\left(1-24\sum^{\infty}_{n=1}\frac{nq^{2n}}{1-q^{2n}}\right)-p\left(1-24\sum^{\infty}_{n=1}\frac{nq^{2p n}}{1-q^{2p n}}\right),
\end{equation}
for $p=3,5,7,9,11,15,17,19,23,25,31,35$. Attached with these evaluations underlines the quite interesting formula 
\begin{equation}
T_{p,r}=\frac{4}{\pi^2}k_r(k'_r)^2K^2\frac{d}{dk_r}\log\left(M_{p}(r)^{-3}\frac{k_rk'_r}{k_{p^2r}k'_{p^2r}}\right),
\end{equation}
where 
\begin{equation}
M_{p}(r)=\frac{K\left(k_{p^2r}\right)}{K\left(k_r\right)}=\frac{K[p^2r]}{K[r]}
\end{equation}
is the multiplier and $k'_r=\sqrt{1-k_r^2}$. The multiplier is also algebraic function when $r$ is positive rational and $p$ integer grater than 1. (Note that $K=K(k_r)=K[r]$ and $E=E(k_r)=E[r]$.)\\

For now on we assume the notation of $T_{p,r}$ as given in (166)\\
\\
\textbf{Proposition 2.}\\
If $A_{p,r}$ is given from
\begin{equation}
A_{p,r}:=\frac{f(-q^2)}{q^{c}f(-q^{2p})}\textrm{, }c=\frac{p-1}{12},
\end{equation}
then
\begin{equation}
T_{p,r}=-\frac{24\sqrt{r}}{\pi A_{p,r}}\cdot\frac{dA_{p,r}}{dr}.
\end{equation}
\\
\textbf{Proof.}\\
Take the logarithm of (169) and use (149) directly. The  result follows easily after differentiation with respect to $r$.\\
\\
\textbf{Proposition 3.}\\
If $q=e^{-\pi\sqrt{r}}$, $r>0$, then
\begin{equation}
P(q^2)=1-24\sum^{\infty}_{n=1}\frac{nq^{2n}}{1-q^{2n}}=\frac{3}{\pi\sqrt{r}}+\left(1+k^2_r-\frac{3\alpha(r)}{\sqrt{r}}\right)\frac{4}{\pi^2}K^2[r].
\end{equation}
\\
\textbf{Proof.}\\
Let $q=e^{-\pi\sqrt{r}}$, $r>0$. Differentiating with respect to $r$ the relation
\begin{equation}
\log\left(f(-q^2)\right)=\sum^{\infty}_{n=1}\log\left(1-q^{2n}\right)\textrm{, }|q|<1,
\end{equation}
and using (see [13] Chapter 21, Miscellaneous examples 10, pg. 488):
\begin{equation}
f(-q^2)^6=\prod^{\infty}_{n=1}\left(1-q^{2n}\right)^6=\frac{2kk'K(k)^3}{\pi^3q^{1/2}},
\end{equation}
we get
\begin{equation}
\frac{1}{6}\frac{d}{dr}\log\left(\frac{2k_rk'_rK^3}{\pi^3q^{1/2}}\right)=-2\sum^{\infty}_{n=1}\frac{nq^{2n-1}}{1-q^{2n}}\frac{dq}{dr}.
\end{equation}
After some calculations we arrive to ($P(q^2)=1-24P^{*}(q^2)$):
\begin{equation}
\frac{1}{24}\left(\frac{\pi}{\sqrt{r}}+12\frac{1}{K}\frac{dK}{dk_r}\frac{dk_r}{dr}+\frac{4}{k_r}\frac{dk_r}{dr}+\frac{4}{k'_r}\frac{dk'_r}{dr}\right)=2q^{-1}P^{*}(q^2)\frac{q\pi}{2\sqrt{r}}
\end{equation}
Using the known relations (see [5] Chapter 17 Entry 9 pg. 120 and [14] Chapter 11 Entry 30 pg. 87-88):
\begin{equation}
\frac{dk_r}{dr}=\frac{-k_r(k'_r)^2K^2}{\pi\sqrt{r}},
\end{equation}
\begin{equation}
\frac{dk_r'}{dr}=\frac{k_r^2k_r'K^2}{\pi\sqrt{r}}
\end{equation}
and ([12] Chapter 1, Section 1.3, pg. 7-11):
\begin{equation}
\frac{dK}{dk_r}=\frac{E}{k_r(k_r')^2}-\frac{K}{k_r},
\end{equation}
we arrive to
\begin{equation}
P^{*}(q^2)=-\frac{1}{24}+\frac{K^2}{6\pi^2}+\frac{K^2k_r^2}{6\pi^2}-\frac{\alpha(r)K^2}{2\pi^2\sqrt{r}}+\frac{1}{8\pi\sqrt{r}}.
\end{equation}
From this along with $P(q^2)=1-24P^{*}(q^2)$, we get the result.\\ 
\\
\textbf{Proposition 4.}\\
Let $r>0$ and $q=e^{-\pi\sqrt{r}}$, then
$$
1-24\sum^{\infty}_{n=1}\frac{nq^n}{1-q^n}
=\frac{6}{\pi\sqrt{r}}+\left(1+k^2_r-\frac{6 \alpha(r)}{\sqrt{r}}\right)\frac{4K[r]^2}{\pi^2}=
$$
\begin{equation}
=\frac{6}{\pi\sqrt{r}}+s_1(r)\theta_3(q)^4,
\end{equation}
where
\begin{equation}
s_1(r):=1-\frac{6\alpha(r)}{\sqrt{r}}+k_r^2
\end{equation}
\\
\textbf{Proof.}\\
From (171) setting $r\rightarrow r/4$ and using
\begin{equation}
k_{r/4}=\frac{2\sqrt{k_r}}{1+k_r}\textrm{, }M_2(r)=\frac{1+k'_r}{2}
\end{equation}
and (see [12]) 
\begin{equation}
\alpha(4r)=(1+k_{4r})^2\alpha(r)-2\sqrt{r}k_{4r}, 
\end{equation}
we get the result.\\
\\
\textbf{Proposition 5.}\\
If $q=e^{-\pi\sqrt{r}}$, $r>0$, then
\begin{equation}
T_{2,r}=-\frac{\theta_3(q)^4}{2}-\frac{\theta_4(q)^4}{2}=-\frac{1}{2}\frac{\eta(q^2)^{20}}{\eta(q)^8\eta(q^4)^8}-\frac{1}{2}\frac{\eta(q)^8}{\eta(q^2)^4}
\end{equation}
\\
\textbf{Proof.}\\
From relations (171),(183) we get
$$
T_{2,r}=1-24\sum^{\infty}_{n=1}\frac{nq^{2n}}{1-q^{2n}}-2\left(1-24\sum^{\infty}_{n=1}\frac{nq^{4n}}{1-q^{4n}}\right)=
$$
$$
=-\frac{2K^2}{\pi^2}\left(1+(k'_r)^2\right)=-\frac{\theta_3(q)^4}{2}-\frac{\theta_4(q)^4}{2}.
$$
But also holds
\begin{equation}
\theta_3(q)=\frac{\eta(q^2)^5}{\eta(q)^2\eta(q^4)^2},
\end{equation}
and
\begin{equation}
\theta_4(q)=\frac{\eta(q)^2}{\eta(q^2)}.
\end{equation}
Using the above equations the result follows.\\

It is easy to see someone that 
\begin{equation}
\delta(q):=P(q)-2P(q^2)=-\frac{4}{\pi^2}(1+k_r^2)K[r]^2=-\theta_2(q)^4-\theta_3(q)^4
\end{equation} 
Hence also
\begin{equation}
T_{2p,r/4}=2T_{p,r}+\delta(q)
\end{equation}
and
\begin{equation}
T_{p,4r}-\frac{1}{2}T_{p,r}=-\frac{1}{2}\left(\delta(q^2)-p\delta(q^{2p})\right)
\end{equation}
Moreover the next doublication formula is trivial:
\begin{equation}
T_{2p,r}=T_{p,r}+p\delta(q^{2p}).
\end{equation}
Although relation (188) is more applicable for finding special values.\\
Set now
\begin{equation}
s_2(r):=1-\frac{3\alpha(r)}{\sqrt{r}}+k_r^2
\end{equation}
\\
\textbf{Proposition 6.}\\
If $q=e^{-\pi\sqrt{r}}$, $r>0$ and $1<p\in\textbf{N}$, then
$$
T_{p,r}=\frac{4K[r]^2}{\pi^2\sqrt{r}}\left[\left(3\alpha(p^2r)-p\sqrt{r}(1+k^2_{p^2r})\right) M_p(r)^2-\left(3\alpha(r)-\sqrt{r}(1+k^2_r)\right)\right]=
$$
\begin{equation}
=s_2(r)\theta_3(q)^4-p\cdot s_2(p^2r)\theta_3(q^p)^4
\end{equation}
and
\begin{equation}
\frac{M^2_p}{\sqrt{r}}\cdot\alpha(p^2r)-\frac{1}{\sqrt{r}}\cdot\alpha(r)=\frac{p  M_p^2}{3}(1+k^2_{p^2r})-\frac{1}{3}(1+k^2_r)+\frac{\pi^2 T_{p,r}}{12K[r]^2}
\end{equation}
\\
\textbf{Proof.}\\
Set $r\rightarrow p^2 r$ in (171) of Proposition 3 and make simplifications to get
\begin{equation}
1-24\sum^{\infty}_{n=1}\frac{nq^{2pn}}{1-q^{2pn}}=
\frac{3}{\pi\sqrt{r}p}+\frac{4K[r]^2}{\pi^2\sqrt{r}pM^2_p }\left[-3\alpha(p^2r)+p\sqrt{r}(1+k^2_{p^2r})\right],
\end{equation}
then evaluate $T_{p,r}$.\\
\\
\textbf{Note 1.}\\
The two functions $P(q)$ and $s_2(r)\theta_3(q)^4$ are not the same. As someone can see easily 
\begin{equation}
P(q^2)=s_2(r)\theta_3(q)^4+\frac{3}{\pi\sqrt{r}}\textrm{, }\forall r>0.
\end{equation}
\\
\textbf{Proposition 7.}\\
If $q=e^{-\pi\sqrt{r}}$, $r>0$ and  $(a,b,c)\stackrel{\theta}{\equiv}(1,1,1),(1,2,2),(1,1,2)$, with $D=-3$, $-4$, $-7$ respectively, then we have the following evaluations
\begin{equation}
\phi(a,b,c;q^2)=\sqrt{\frac{T_{|D|,r}}{1-|D|}},
\end{equation}
where $T_{|D|,r}=T_{p,r}$ is given from (45). More precisely the cases $p=3,4,7$ are solvable.\\   
\\
\textbf{Proposition 8.}\\
If $a_1,b_1,c_1$ are given integers with $a_1>0$, $c_1>0$, $D=b_1^2-4a_1c_1<0$ and exist $a_2,c_2$ integers with $a_2>0$, $c_2>0$ such  $(a_1,b_1,c_1)\stackrel{\theta}{\equiv}(a_2,0,c_2)$, then if $q=e(z)$, $Im(z)>0$, we have $\phi(a_1,b_1,c_1;q)^2\in M_2\left(\Gamma_0(N)\right)$ and $dim\left(M_2\left(\Gamma_0(N)\right)\right)=d=d_1+d_2$ and exist constants $C^{(1)}_l$, $C^{(2)}_l$ such that 
\begin{equation}
\theta_3\left(q^{a_2}\right)^2\theta_3\left(q^{c_2}\right)^2=\sum_{l|N}C^{(1)}_lt_{l}(q)+\sum^{d_2}_{l=1}C^{(2)}_lf_l(q),
\end{equation}
where $f_l(q)$, $l=1,2,\ldots,d_2$ are the elements of the subspace of cusp forms $S_2\left(\Gamma_0(N)\right)$ of $M_2\left(\Gamma_0(N)\right)$ and   $t_{l}(q)=P\left(q\right)-lP\left(q^{l}\right)=E_2(z)-lE_2\left(lz\right)$. In case $z=i\sqrt{r}$, $r>0$, $q=e^{-2\pi\sqrt{r}}$, then $t_{l}(q)=T_{l,r}$ and we can write
\begin{equation}
\theta_3\left(q^{2a_2}\right)^2\theta_3\left(q^{2c_2}\right)^2=\sum_{l|N}C^{(1)}_lT_{l,r}+\sum^{d_2}_{l=1}C^{(2)}_lf_l(q^2)\textrm{, }q=e^{-\pi\sqrt{r}}\textrm{, }r>0.
\end{equation}
Relation (198) can also be used to evaluate $\alpha(|D|^2r)$.\\
\\
\textbf{Application 1.}\\
Assume that $D=-4$, then $Q(x,y)=x^2+y^2$ and $Q$ have array $A=\left(
\begin{array}{cc}
	2\textrm{ }0\\
	0\textrm{ }2
\end{array}
\right)$. The level $N$ is 4 and  $dim\left(M_2\left(\Gamma_0(4)\right)\right)=dim\left(E_2\left(\Gamma_0(4)\right)\right)=1$, hence we can write $\phi(1,0,1;q)^2=C_1\left(E_2(z)-4E_2\left(4z\right)\right)$, (here $C_1=-1/3$) and $q=e(z)$, $Im(z)>0$. Hence for $q=e^{-\pi\sqrt{r}}$, $r>0$
\begin{equation}
\phi\left(1,0,1;q^2\right)^2=-\frac{1}{3}T_{4,r}.
\end{equation}
We use Proposition 8 with
$a_1=c_1=1$, $b_1=0$ and $a_2=1$, $b_2=c_2=0$ to write
\begin{equation}
\phi\left(1,0,1;q^2\right)^2=\theta_3\left(q^2\right)^4=\frac{4K[4r]^2}{\pi^2}=\left(1+k'_r\right)^2\frac{K^2}{\pi^2}.
\end{equation}
Hence using relation (193) of Proposition 6 with $p=4$ we get:
$$
\frac{M^2_4}{\sqrt{r}}\cdot\alpha(16r)-\frac{1}{\sqrt{r}}\cdot\alpha(r)=\frac{4  M_4^2}{3}(1+k^2_{16r})-\frac{1}{3}(1+k^2_r)+\frac{\pi^2 T_{4,r}}{12K[r]^2}.
$$
Or equivalently
\begin{equation}
\frac{M^2_4}{\sqrt{r}}\cdot\alpha(16r)-\frac{1}{\sqrt{r}}\cdot\alpha(r)=\frac{4  M_4^2}{3}(1+k^2_{16r})-\frac{1}{3}(1+k^2_r)-\frac{(1+k'_r)^2}{4},
\end{equation}
where
\begin{equation}
M_4(r)=\frac{1+k'_r}{2}\frac{1+k'_{4r}}{2}=\frac{(1+\sqrt{k'_r})^2}{4},
\end{equation}
since 
$$
k_{4r}=\frac{1-k'_r}{1+k'_r}.
$$
After simplifications we get the following:\\
\\
\textbf{Proposition 9.}\\
If $r>0$, then
\begin{equation}
\alpha(16r)=\frac{16\alpha(r)-4(1-k'_r)(3+k'_r)\sqrt{r}}{(1+\sqrt{k'_r})^4}
\end{equation}
\\
and\\
\\
\textbf{Proposition 10.}\\
We have $q=e^{-\pi\sqrt{r}}$, $r>0$
\begin{equation}
T_{4,r}=-3\theta_3(q^2)^4=-3\frac{\eta(q^4)^{20}}{\eta(q^2)^8\eta(q^8)^8}.
\end{equation}
\\
\textbf{Proof.}\\ 
The second equality of (204) follows from (185).\\
\\
\textbf{Application 2.}\\
Now we assume the theta function
\begin{equation}
\phi\left(1,0,2;q\right)^2=\left(\sum^{\infty}_{n,m=-\infty}q^{n^2+2m^2}\right)^2.
\end{equation}
The level of the function is $N=8$. But in $M_2\left(\Gamma_0(8)\right)$ we have $dimE_2\left(\Gamma_0(8)\right)=3$ and $dimS_2\left(\Gamma_0(8)\right)=1$. Using the method of Lemma 2, pg.14-15, we find that the cusp form is
\begin{equation}
\Delta_8(q)=\frac{\eta\left(q^4\right)^8}{\eta\left(q^2\right)^4}.
\end{equation} 
Also the expansion of (205) into modular forms of weight 2 in $M_2\left(\Gamma_0\left(8\right)\right)$ is
\begin{equation}
\phi(1,0,2;q)^2=\frac{1}{18}\left(E_2(q)-4E_2\left(q^4\right)\right)-\frac{1}{6}\left(E_2(q)-8E_2\left(q^8\right)\right)+\frac{4}{3}\Delta_8(q).
\end{equation}
Or setting $q\rightarrow q^2$, $q=e^{-\pi\sqrt{r}}$, $r>0$:
$$
\theta_3\left(q^2\right)^2\theta_3\left(q^4\right)^2=\frac{1}{18}T_{4,r}-\frac{1}{6}T_{8,r}+\frac{4}{3}\Delta_8(q^2).
$$
Using relation (204) we get
\begin{equation}
T_{8,r}=-\theta_3\left(q^2\right)^4-6\theta_3\left(q^2\right)^2\theta_3\left(q^4\right)^2+8 \Delta_8\left(q^2\right).
\end{equation}   
Hence we have the next\\
\\
\textbf{Proposition 11.}\\
If $q=e^{-\pi\sqrt{r}}$, $r>0$, then
\begin{equation}
T_{8,r}=\frac{8 \eta\left(q^8\right)^8}{\eta\left(q^4\right)^4}-\frac{\eta\left(q^4\right)^{20}}{\eta\left(q^2\right)^8 \eta\left(q^8\right)^8}-\frac{6\eta\left(q^4\right)^6 \eta\left(q^8\right)^6}{\eta\left(q^2\right)^4 \eta\left(q^{16}\right)^4}.
\end{equation}
\\
\textbf{Application 3.}\\
Assume $D=-20$, then one choice is $Q(x,y)=x^2+2xy+6y^2$ and $Q$ have array $A=\left(
\begin{array}{cc}
	2\textrm{ }\textrm{ }2\\
	2\textrm{ }12
\end{array}
\right)$. The level $N$ is 20 and  $dim\left(E_2\left(\Gamma_0(20)\right)\right)=5$, $dim\left(S_2\left(\Gamma_0(20)\right)\right)=1$, hence we can write (having in mind relation (73)):
$$
\phi(1,2,6;q)^2=\theta_3(q)^2\theta_3(q^5)^2=-\frac{1}{20}(E_2(q)-2E_2(q^2))+\frac{1}{6}(E_2(q)-5E_2(q^5))-
$$
\begin{equation}
-\frac{1}{6}(E_2(q)-10E_2(q^{10}))-\frac{1}{6}\Delta_{20}(q),
\end{equation}
where
\begin{equation}
\Delta_{20}(q)=\frac{\eta(q^4)^4\eta(q^{20})^4}{\eta(q^2)^2\eta(q^{10})^2}.
\end{equation}
Hence
$$
\theta_3(q^2)^2\theta_3(q^{10})^2=-\frac{1}{6}T_{2,r}+\frac{1}{6}T_{5,r}-\frac{1}{6}T_{10,r}-16\Delta_{20}(q^2)
$$
Or using (204) and 
\begin{equation}
T_{5,r}=-4\eta(q^2)^4A_{5,r}^{-5}\sqrt{125+22 A_{5,r}^{6}+A_{5,r}^{12}},
\end{equation}
where 
\begin{equation}
A_{5,r}=\frac{\eta(q^2)}{\eta(q^{10})}\textrm{, }q^{-\pi\sqrt{r}}\textrm{, }r>0, 
\end{equation}
we get:\\
\\
\textbf{Proposition 12.}\\
If $q=e^{-\pi\sqrt{r}}$, $r>0$, then $T_{5,r}<0$ and
\begin{equation}
T^2_{5,r}=T^2_5(q^2)=16\frac{\eta(q^2)^{12}+22\eta(q^2)^6\eta(q^{10})^6+125\eta(q^{10})^{12}}{\eta(q^2)^2\eta(q^{10})^2}
\end{equation}
$$
T_{10}(q^2)=T_{5}(q^2)+\frac{\eta(q)^8}{2 \eta\left(q^2\right)^4}+\frac{\eta\left(q^2\right)^{20}}{2 \eta(q)^8 \eta\left(q^4\right)^8}
-\frac{96 \eta\left(q^8\right)^4 \eta\left(q^{40}\right)^4}{\eta\left(q^4\right)^2 \eta\left(q^{20}\right)^2}-
$$
\begin{equation}
-\frac{6 \eta\left(q^4\right)^{10} \eta\left(q^{20}\right)^{10}}{\eta\left(q^2\right)^4 \eta\left(q^8\right)^4 \eta\left(q^{10}\right)^4 \eta\left(q^{40}\right)^4}.
\end{equation}
Also
\begin{equation}
\theta_3(q^2)^2\theta_3(q^6)^2=-\frac{1}{6}T_{2,r}-\frac{1}{6}T_{3,r}+\frac{1}{6}T_{4,r}+\frac{1}{6}T_{6,r}-\frac{1}{6}T_{12,r}.
\end{equation}
\\
\textbf{Proposition 13.}\\
If $q=e^{-\pi\sqrt{r}}$, $r>0$, then
\begin{equation}
T_{3,r}=-2\left(\frac{3\eta(q)^2\eta(q^4)^2\eta(q^6)^{15}}{4\eta(q^2)^5\eta(q^3)^6\eta(q^{12})}+\frac{\eta(q^2)^{15}\eta(q^3)^2\eta(q^{12})^2}{4\eta(q)^6\eta(q^4)^6\eta(q^6)^5}\right)^2.
\end{equation}
Also if $|q|<1$, then
\begin{equation}
\sum^{\infty}_{n,m=-\infty}q^{2n^2+2nm+2m^2}=\frac{3\eta(q)^2\eta(q^4)^2\eta(q^6)^{15}}{4\eta(q^2)^5\eta(q^3)^6\eta(q^{12})}+\frac{\eta(q^2)^{15}\eta(q^3)^2\eta(q^{12})^2}{4\eta(q)^6\eta(q^4)^6\eta(q^6)^5}
\end{equation}
and
\begin{equation}
\sum^{\infty}_{n,m=-\infty}q^{2n^2+2nm+2m^2}=\frac{\theta_3(q)^4+3\theta_3(q^3)^4}{4\theta_3(q)\theta_3(q^3)}.
\end{equation}
\\
\textbf{Proof.}\\
For $|q|<1$, we have
\begin{equation}
\theta_3(q)=\frac{\eta(q^2)^5}{\eta(q)^2\eta(q^4)^2}.
\end{equation}
Also (see [5] last chapter):
\begin{equation}
\left(\sum^{\infty}_{n,m=-\infty}q^{2n^2+2nm+2m^2}\right)^2=-\frac{1}{2}T_{3,r}=\left(\frac{\theta_3(q)^4+3\theta_3(q^3)^4}{4\theta_3(q)\theta_3(q^3)}\right)^2.
\end{equation}
Combining the above relations we get the two results.\\
\\
\textbf{Proposition 14.}\\
If $q=e^{-\pi\sqrt{r}}$, $r>0$, then
\begin{equation}
-\frac{1}{6}T_{7,r}=\left(\frac{\eta(q^4)^5\eta(q^{28})^5}{\eta(q^2)^2\eta(q^8)^2\eta(q^{14})^2\eta(q^{56})^2}+4\frac{\eta(q^8)^2\eta(q^{56})^2}{\eta(q^4)\eta(q^{28})}\right)^2.
\end{equation}
If $|q|<1$, then
\begin{equation}
\sum^{\infty}_{n,m=-\infty}q^{n^2+nm+2m^2}=\frac{\eta(q^2)^5\eta(q^{14})^5}{\eta(q)^2\eta(q^4)^2\eta(q^{7})^2\eta(q^{28})^2}+4\frac{\eta(q^4)^2\eta(q^{28})^2}{\eta(q^2)\eta(q^{14})}.
\end{equation}
and
\begin{equation}
\sum^{\infty}_{n,m=-\infty}q^{n^2+nm+2m^2}=\theta_3(q)\theta_3(q^{7})+4q^2\psi(q^2)\psi(q^{14}),
\end{equation}
where $\psi(q):=\sum^{\infty}_{n=0}q^{n(n+1)/2}$, $|q|<1$.\\
\\
\textbf{Proof.}\\
Ramanujan have proved that for $|q|<1$ holds (see [5] last chapter): 
\begin{equation}
-\frac{1}{6}T_{7}(q^2)=\left(\theta_3(q^2)\theta_3(q^{14})+4q^4\psi(q^4)\psi(q^{28})\right)^2.
\end{equation} 
But also it is
\begin{equation}
\psi(q)=\frac{f(-q^2)^2}{f(-q)}=q^{-1/8}\frac{\eta(q^2)^2}{\eta(q)}.
\end{equation}
Hence from (225),(226),(220) we get (222). Also from (222) and (78) we get (223) and (224).\\
\\
\textbf{Proposition 15.}\\
If $q=e^{-\pi\sqrt{r}}$, $r>0$, then
$$
-\frac{1}{2}T_{11,r}=-\frac{1}{2}T_{11}(q^2)=20\eta(-q)^2\eta(-q^{11})^2+32\eta(q^2)^2\eta(q^{22})^2+
$$
\begin{equation}
+20\frac{\eta(q^2)^4\eta(q^{22})^4}{\eta(-q)^2\eta(-q^{11})^2}
+5\frac{\eta(q^2)^{10}\eta(q^{22})^{10}}{\eta(q)^4\eta(q^4)^4\eta(q^{11})^4\eta(q^{44})^4}
\end{equation}
and
\begin{equation}
\left(\sum^{\infty}_{n,m=-\infty}q^{n^2+n m+3m^2}\right)^2=-\frac{1}{10}T_{11}(q)+\frac{8}{5}\eta(q)^2\eta(q^{11})^2.
\end{equation}
\\
\textbf{Proof.}\\
Ramanujan have proved that for $|q|<1$ holds (see [5] last chapter):
$$
-\frac{1}{2}T_{11}(q^2)=5\theta_3(q)^2\theta_3(q^{11})^2-20qf(q)^2f(q^{11})^2
+32q^2f(-q^2)^2f(-q^{22})^2-
$$
\begin{equation}
-20q^3\psi(-q)^2\psi(-q^{11})^2.
\end{equation}
Using the above equation and (226),(220), we get (227). Also from (88) we get (228).\\ 
\\
\textbf{Proposition 16.}\\
The functions $A_{p,r}$ and $T_{p,r}\theta_3\left(q\right)^{-4}$ are algebraic and $\theta_3(q)^p$, $q=e^{-\pi\sqrt{r}}$ transcendental when $r,p\in\textbf{Q}^{*}_{+}$.\\ 
\\
\textbf{Proof.}\\
From (169),(173) we get
$$
\left(A_{p,r}\right)^6=\frac{f\left(-q^2\right)^6}{q^{\frac{p-1}{2}}f\left(-q^{2p}\right)^6}=\frac{2k_rk'_rK[r]^3}{\pi^3q^{1/2}q^{\frac{p-1}{2}}2k_{p^2r}k'_{p^2r}K[p^2r]^3\frac{1}{\pi^3 q^{p/2}}}.
$$
Hence using (168)
\begin{equation}
A_{p,r}=M_{p}(r)^{-1/2}\sqrt[6]{\frac{k_rk'_r}{k_{p^2r}k'_{p^2r}}}
\end{equation}
and the first result follows.\\
The algebricity of $T_{p,r}\theta_3\left(q\right)^{-4}$ follows from the expansion (192) for $\frac{\pi^2}{K[r]^2} T_{p,r}$ and relation
\begin{equation}
\theta_3(q)^2=\frac{2K[r]}{\pi}.
\end{equation}  
The transcendence of $\theta_3\left(e^{-\pi\sqrt{r}}\right)$ follows from the fact that of $K[r]$ is transcendental and algebraically independent from $\pi$ when $r$ is positive rational (see [5]).

\section{Evaluating the quotients $A_{p,r}$ for $p=3,7$ using known integral functions}

In this section we give evaluations of integrals and hypergeometric functions using the quotients $A_{p,r}$ and Carty's function $\Pi(r)$ for the cases $p=3,7$.\\ 
\\
\textbf{Theorem 22.}\\
If $r>0$, then
\begin{equation}
\frac{6}{A_{3,r}^2}\cdot {}_2F_{1}\left(\frac{1}{6},\frac{2}{3};\frac{7}{6};-\frac{27}{A_{3,r}^{12}}\right)=3\sqrt[3]{2k_{4r}}\cdot {}_2F_1\left(\frac{1}{3},\frac{1}{6};\frac{7}{6};k_{4r}^2\right)=\Pi(r)
\end{equation}
\\
\textbf{Proof.}\\
From [5] chapter 21 Entry 3 (i) and (170) we have
$$
\frac{dA_{3,r}}{dr}=\frac{\pi}{12\sqrt{r}} q^{1/3}f\left(-q^2\right)^4A^{-5}_{3,r}\left(27+A_{3,r}^{12}\right)^{2/3}=
$$
\begin{equation}
=\frac{\pi}{12\sqrt{r}}\eta_1\left(i\sqrt{r}\right)^4 A^{-5}_{3,r}\left(27+A_{3,r}^{12}\right)^{2/3}.
\end{equation}
Integrating the above equation we get:
\begin{equation}
\int^{+\infty}_{A_{3,r}}\frac{t^5}{(27+t^{12})^{2/3}}dt=\frac{\pi}{6}\int^{+\infty}_{\sqrt{r}}\eta_1(it)^4dt
\end{equation}
We rewrite the above as
$$
\int^{+\infty}_{A_{3,r/4}}\frac{t^5}{(27+t^{12})^{2/3}}dt=\frac{\pi}{12}\int^{+\infty}_{\sqrt{r}}\eta_1(it/2)^4dt
$$
But it is known that the function $m(x)$ defined as
\begin{equation}
x=\pi\int^{+\infty}_{\sqrt{m(x)}}\eta_1\left(it/2\right)^4dt
\end{equation}
has inverse 
\begin{equation}
m^{(-1)}(r)=b_r=3\sqrt[3]{2k_r}\cdot {}_2F_1\left(\frac{1}{3},\frac{1}{6};\frac{7}{6};k_r^2\right). 
\end{equation}
Hence
$$
\int^{+\infty}_{A_{3,m(x)/4}}\frac{t^5}{(27+t^{12})^{2/3}}dt=\frac{\pi}{12}\int^{+\infty}_{\sqrt{m(x)}}\eta(it/2)^4dt=\frac{x}{12}
$$
and consequently
$$
\int^{+\infty}_{A_{3,r/4}}\frac{t^5}{(27+t^{12})^{2/3}}dt=\frac{1}{4}\sqrt[3]{2k_r}\cdot {}_2F_1\left(\frac{1}{3},\frac{1}{6};\frac{7}{6};k_r^2\right).
$$
Having in mind that
$$
\int^{+\infty}_{X}\frac{t^5}{(27+t^{12})^{2/3}}dt=\frac{1}{2X^2}\cdot {}_2F_{1}\left(\frac{1}{6},\frac{2}{3};\frac{7}{6};-\frac{27}{X^{12}}\right),
$$  
we get the result.\\
\\
\textbf{Theorem 23.}\\
Set $q=e^{-\pi\sqrt{r}}$, $r>0$, then if
\begin{equation}
A_{7,r}=\frac{f(-q^2)}{q^{1/2}f(-q^{14})},
\end{equation}
we have
\begin{equation}
\int^{+\infty}_{A_{7,r}^4}\frac{t^{1/6}}{(49+13t+t^2)^{2/3}}dt=\Pi(r)
\end{equation}
\\
\textbf{Proof.}\\
Again from [5] chapter (21) Entry 5 (i) we have
\begin{equation}
\frac{dA_{7,r}}{dr}=\frac{\pi}{4\sqrt{r}} q^{1/3}f\left(-q^2\right)^4A^{-11/3}_{7,r}\left(49+13 A_{7,r}^4+A_{7,r}^8\right)^{2/3}
\end{equation}
\begin{equation}
\int^{+\infty}_{A_{7,r}}\frac{t^{11/3}}{4\left(49+13 t^4+t^8\right)^{2/3}}dt=\frac{3}{4}\sqrt[3]{2k_{4r}}\cdot{ }_2F_1\left(\frac{1}{3},\frac{1}{6};\frac{7}{6};k_{4r}^2\right).
\end{equation}

Assume that $q=e^{-\pi\sqrt{r}}$, $r>0$, then we re-set
\begin{equation}
x_p=x_p(r)=A_{p,r}=\frac{\eta_1\left(i\sqrt{r}\right)}{\eta_1\left(i p\sqrt{r}\right)}.
\end{equation}
and
\begin{equation}
z_p=z_p(r)=T_{p,r}=P(q)-pP\left(q^p\right)
\end{equation}
The function $F_7(x)$ is such that
\begin{equation}
\int^{+\infty}_{F_7^4(x)}\frac{t^{1/6}}{(49+13t+t^2)^{2/3}}dt=x
\end{equation}
and
\begin{equation}
x_7=A_{7,r}=F_7\left(\Pi(r)\right).
\end{equation}
From (170) we have in general\\
\\
\textbf{Theorem 24.}\\
If $q=e^{-\pi\sqrt{r}}$, $r>0$, then exists a function $F_p(x)$ such that
\begin{equation}
-\frac{\pi}{24\sqrt{r}}\left(P(q^2)-pP\left(q^{2p}\right)\right)=\frac{d}{dr}\log F_p\left(\Pi\left(r\right)\right).
\end{equation}

For $p=7$ we get integrating\\
\\
\textbf{Theorem 25.}\\
If $q=e^{-\pi t}$, $t>0$, then
\begin{equation}
\frac{\pi}{4} \int_{2\sqrt{r}}^{+\infty}\left(\sum^{\infty}_{n,m=-\infty}q^{n^2+nm+2 m^2}\right)^2dt=\log\left(F_7\left(\Pi\left(r\right)\right)\right)=\log x_7.
\end{equation}

Continuing our arguments Proposition 2 give us
\begin{equation}
\frac{-\pi}{24\sqrt{r}}z_p(r)=\frac{1}{x_p(r)}\frac{d}{dr}x_p(r).
\end{equation}
Also if 
\begin{equation}
A_{p,r}=x_p(r)=F_p\left(\Pi(r)\right)
\end{equation}
and
\begin{equation}
T_{p,r}=z_p(r)=Q_{p}\left(x_p(r)\right),
\end{equation}
then
\begin{equation}
\frac{-24\sqrt{r}}{\pi}\frac{1}{x_p(r)}\frac{dx_p(r)}{dr}=Q_p\left(x_p(r)\right)
\end{equation}
and inverting this last equation we obtain
\begin{equation}
\frac{-24\sqrt{x_p^{(-1)}(t)}}{\pi}\frac{1}{t\cdot x^{(-1)}_p{'}(t)}=Q_p(t)
\end{equation}
Inverting relation (248) we have
\begin{equation}
x_p^{(-1)}(t)=m\left(F_p^{(-1)}(t)\right).
\end{equation}
But $m\left(\Pi(t)\right)=t$. Hence differentiating we get $
m'\left(\Pi(t)\right)\Pi'\left(t\right)=1$. Hence
\begin{equation} m'(t)=-\frac{\sqrt{m(t)}}{\pi\eta_1\left(i\sqrt{m(t)}\right)^4},
\end{equation}
since it holds
\begin{equation}
\frac{d\Pi(r)}{dr}=-\pi\frac{\eta_1\left(i\sqrt{r}\right)^4}{\sqrt{r}}.
\end{equation}
Hence
\begin{equation}
Q_p(t)=\frac{24}{t\cdot F^{(-1)}_p{'}(t)}\eta_1\left(i\sqrt{m\left(F^{(-1)}_p(t)\right)}\right)^4
\end{equation} 
and we get the next\\
\\
\textbf{Theorem 26.}\\
\begin{equation}
z_p=\frac{24}{x_p\cdot F^{(-1)}_p{'}(x_p)}\eta_1\left(i\sqrt{m\left(F^{(-1)}_p(x_p)\right)}\right)^4
\end{equation}
where
\begin{equation}
x_p=\frac{\eta_1\left(i\sqrt{r}\right)}{\eta_1\left(i p \sqrt{r}\right)}\textrm{, }r>0.
\end{equation}
In cases $p=3,5,7$, we have
\begin{equation}
F_3^{(-1)}(x)=\frac{6}{x^2}\cdot{ }_{2}F_1\left(\frac{1}{6},\frac{2}{3};\frac{7}{6};\frac{27}{x^{12}}\right),
\end{equation}
\begin{equation}
F_5^{(-1)}(x)=\int^{+\infty}_{x}\frac{dt}{t^{1/4}\sqrt{125+22t+t^2}}
\end{equation}
and
\begin{equation}
F_7^{(-1)}(x)=\int^{+\infty}_{x^4}\frac{t^{1/6}}{(49+13t+t^2)^{2/3}}dt.
\end{equation}

\newpage

\[
\]

\centerline{\bf References}\vskip .2in

[1]: J.V. Armitage, W.F. Eberlein. 'Elliptic Functions'. Cambridge University Press. (2006)

[2]: N.D. Bagis, M.L. Glasser. 'Conjectures on the Evaluation of Alternative Modular Bases and Formulas Approximating $1/\pi$'. Journal of Number Theory. (Elsevier), (2012).

[3]: N.D. Bagis, M.L. Glasser. 'Conjectures on the evaluation of certain functions with algebraic properties'. Journal of Number Theory. 155 (2015), 63-84

[4]: N.D. Bagis. 'On the Complete Evaluation of Jacobi Theta Functions'. arXiv:1503.01141v1 [math.GM],(2015).

[5]: B.C. Berndt. 'Ramanujan`s Notebooks Part III'. Springer Verlag, New York (1991).

[6]: B.C. Berndt, S. Bhargava, F.G. Garvan. 'Ramanujan's theories of elliptic functions to alternative bases'. Trans. Amer. Math. Soc. 347 (1995) 4163-4244.

[7]: J.M. Borwein, M.L. Glasser, R.C. McPhedran, J.G. Wan, I.J. Zucker. 'Lattice Sums Then and Now'. Cambridge University Press. New York, (2013).

[8]: Shaun Cooper, Dongxi Ye. 'Level 14 and 15 analogues of Ramanujan's elliptic functions to alternative bases'. Transactions of the American Mathematical Society, \textbf{368} (2016), 7883-7910.

[9]: Don. Zagier. 'Elliptic Modular Forms and Their Applications'. Available from zagier@mpim-bonn.mpg.de

[10]: Toshijune Miyake. 'Modular Forms'. Springer Verlang, Berlin, Heidelberg, (1989).

[11]: Carlos J. Moreno, Samuel S. Wagstaff. Jr. 'Sums of Quares of Integers'. Chapman and Hall/CRC, Taylor and Francis Group, (2006).

[12] J.M. Borwein and P.B. Borwein. 'Pi and the AGM: A Study in Analytic Number Theory and Computational Complexity', Wiley, New York, 1987.

[13] E.T. Whittaker and G.N. Watson. 'A course on Modern Analysis'. Cambridge U.P. 1927.

[14] Bruce. C. Berndt. 'Ramanujan`s Notebooks Part II'. Springer-Verlag, New York. 1989.

[15] N.D. Bagis and M.L. Glasser. 'On the Transcendence of Complete Elliptic Integrals of the First Kind and Values of the Gamma Function'. submitted

\end{document}